\documentclass[a4paper]{article}
\usepackage{typearea}
\usepackage{fix-cm}
\usepackage[noblocks]{authblk}

\usepackage{graphicx}
\usepackage{amsmath,amsfonts,amssymb}
\usepackage{mathtools}
\usepackage{fixmath}
\usepackage{latexsym}
\usepackage{hyperref}
\usepackage{subcaption}
\usepackage{pgfplots}
\pgfplotsset{compat=1.14,
	     width=5.9cm,
	     height=5.9cm
	    }
\usepgfplotslibrary{groupplots}
\usepackage[export]{adjustbox}

\newtheorem{remark}{Remark}

\newcommand{\uext}{u_{\text{ext}}}
\newcommand{\uinit}{u_{\text{init}}}
\newcommand{\uh}[1]{u_{h}^{#1}}

\newcommand{\urmu}[1]{u_{r,\mu}^{#1}}
\newcommand{\huh}[1]{\hat{u}_{h}^{#1}}
\newcommand{\huhmu}[1]{\hat{u}_{h,\mu}^{#1}}
\newcommand{\hurmu}[1]{\hat{u}_{r,\mu}^{#1}}
\newcommand{\vh}{v_{h}}
\newcommand{\hvh}{\hat{v}_{h}}
\newcommand{\hvr}{\hat{v}_{r}}

\newcommand{\F}{\mathbf{F}}
\renewcommand{\P}{\mathbf{P}}
\newcommand{\I}{\mathbf{I}}

\newcommand{\bs}{\mathbf{s}}
\newcommand{\bx}{\mathbf{x}}

\newcommand{\bPsi}{\mathbold{\Psi}}

\newcommand{\hbpsi}{\hat{\mathbold{\psi}}}
\newcommand{\hbeta}{\hat{\mathbold{\eta}}}
\newcommand{\Psiinit}{\mathbold{\Psi}_{\text{init}}}
\newcommand{\Psih}[1]{\mathbold{\Psi}_{h}^{#1}}
\newcommand{\Psihmu}[1]{\mathbold{\Psi}_{h,\mu}^{#1}}
\newcommand{\Psirmu}[1]{\mathbold{\Psi}_{r,\mu}^{#1}}

\newcommand{\velsymb}{w}
\newcommand{\bvelsymb}{\mathbf{\velsymb}}
\newcommand{\n}{\mathbf{n}}

\newcommand{\Uh}{\hat{U}_h}
\newcommand{\bUh}{\mathbf{\hat{U}}_h}
\newcommand{\bUhG}{\mathbf{\hat{U}}_h^{\Gamma}}

\newcommand{\bM}{\mathbf{M}}

\newcommand{\bUrPsi}{\mathbf{\hat{U}}_{r}}
\newcommand{\bPrPsi}{\mathbf{P}_{r}}
\newcommand{\bUrG}{\mathbf{\hat{U}}_{r}^\Gamma}

\newcommand{\Uru}{\hat{U}_{r}}
\newcommand{\Pru}{P_{r}}

\newcommand{\q}{\bvelsymb}
\newcommand{\qh}[1]{\bvelsymb_{h}^{#1}}
\newcommand{\refqh}[1]{\hat{\bvelsymb}_{h}^{#1}}
\newcommand{\refqhmu}[1]{\hat{\bvelsymb}_{h,\mu}^{#1}}
\newcommand{\refqrmu}[1]{\hat{\bvelsymb}_{r,\mu}^{#1}}
\newcommand{\hn}{\hat{\n}}
\newcommand{\qb}{{\velsymb_\Gamma}}
\newcommand{\qbh}[1]{\bvelsymb_{\Gamma,h}^{#1}}
\newcommand{\refqbh}[1]{\hat{\bvelsymb}_{\Gamma,h}^{#1}}
\newcommand{\refqbhmu}[1]{\hat{\bvelsymb}_{\Gamma,h,\mu}^{#1}}
\newcommand{\refqbrmu}[1]{\hat{\bvelsymb}_{\Gamma,r,\mu}^{#1}}

\newcommand{\T}{\mathcal{T}}
\newcommand{\id}{\mathbf{i}\mathbf{d}}
\newcommand{\refuh}[1]{\hat{u}_h^{#1}}

\newcommand{\mqa}{a_1}
\newcommand{\cmqa}{c_1\!}

\newcommand{\mqb}{a_2}
\newcommand{\cmqb}{{c}_2\!}

\newcommand{\muu}{a_3}
\newcommand{\cmuu}{c_3\!}

\newcommand{\aua}{a_5}
\newcommand{\caua}{{c}_5\!}

\newcommand{\aub}{a_4}
\newcommand{\caub}{{c}_4\!}

\newcommand{\lqa}{l_1}
\newcommand{\clqa}{{c}_6\!}

\newcommand{\lqc}{l_2}
\newcommand{\clqc}{{c}_7}

\newcommand{\hext}{\mathcal{E}}
\newcommand{\hexth}{\mathcal{E}_h}
\begin{document}

\title{Mass Conservative Reduced Order Modeling of a Free Boundary Osmotic Cell Swelling Problem}

\author[a]{Christoph Lehrenfeld}
\author[b]{Stephan Rave}

\affil[a]{Institute for Numerical and Applied Mathematics, University of G\"ottingen, lehrenfeld@math.uni-goettingen.de}
\affil[b]{Applied Mathematics, University of M\"unster, stephan.rave@uni-muenster.de}

\date{May 4, 2018}

\maketitle

\begin{abstract}
 We consider model order reduction for a free boundary problem of an osmotic cell that is
 parameterized by material parameters as well as the initial shape of the cell.
 Our approach is based on an Arbitrary-Lagrangian-Eulerian description of the
 model that is discretized by a mass-conservative finite element scheme.
 Using reduced basis techniques and empirical interpolation, we construct a parameterized reduced
 order model in which the mass conservation property of the full-order model 
 is exactly preserved.
 Numerical experiments are provided that highlight the performance of the resulting reduced order model.
\end{abstract}

\section{Introduction}
Free boundary problems are PDE problems that involve an a priori unknown (free) interface or boundary. These type of problems arise in different applications from physics, engineering, finance and biology. Let us mention a few important application fields where free boundary problems play an important role. In physics and engineering, many situations where different fluids (or solids) are involved, e.g. water and oil in petroleum problems, can be cast into free boundary problems as in the classical Stefan problem \cite{stefan1891theorie}; In finance, optimal stopping of stochastic processes is often solved by reduction to free boundary problems \cite{peskir2006optimal}, and in biology, the mathematical modeling of problems like tumor growth and wound healing leads to free boundary problems \cite{friedman2015free}.

While projection-based reduced order modeling techniques such as reduced basis methods have been successfully
applied to a wide variety of PDE models (see, e.g., \cite{QuarteroniManzoniEtAl2016,HesthavenRozzaEtAl2016,BennerOhlbergerEtAl2017}
for an overview), the reduction of problems that involve an evolving geometry $\Omega(t)$ remains challenging:
Whereas traditional reduction methods are built around the idea of finding a joint linear approximation space for
the entire manifold of solution state vectors, the solutions $u(t)$ of free boundary problems naturally lie in time-dependent
function spaces $V(t)$ that depend on the a priori unknown evolution of $\Omega(t)$.

Taking an Eulerian point of view, a naive approch to resolve this issue is to consider linear
embeddings $\Lambda(t)$ of $V(t)$ into a larger space $V$ of (discontinuous) functions on some $\Omega^* \supset \bigcup_t \Omega(t)$ by
extending the functions with zero on the complement of $\Omega(t)$, and then search for a reduced approximation space $V_N \subset
V$ for the solution manifold $\mathcal{M}_V := \{\Lambda(t)(u(t))\}$.
However, as is well-known from hyperbolic problems with traveling shocks, the moving jump at the boundary of $\Omega(t)$
leads to a slow decay of the Kolmogorov $n$-widths of $\mathcal{M}_V$ (e.g.\ \cite{OhlbergerRave2016}).
Thus, a good low-dimensional linear approximation space $V_N$ for $\mathcal{M}_V$ cannot exist.

Following ideas from \cite{OhlbergerRave2013}, the aim of this paper is to approach this problem by considering nonlinear
approximations of $\mathcal{M}_V$ where we allow transformations of functions $f$ in $V_N$ of the form $\bPsi.f(x) := f(\bPsi(x))$
with some diffeomorphism $\bPsi$ of $\Omega^*$.
If the functions $\bPsi(t)$ are chosen such that the supports of functions in $\hat{\mathcal{M}}_V := \{\bPsi(t).\Lambda(t)(u(t))\}$ have a fixed
boundary, the $n$-widths of $\hat{\mathcal{M}}_V$ will decay fast. Given a good approximation space $V_N$ of
$\hat{\mathcal{M}}_V$, $\bPsi(t)^{-1}.V_N$ will then yield good approximations of $\Lambda(t)(u(t))$.
Noting that such transformations $\bPsi(t)$ induce mappings of a fixed reference domain $\hat\Omega$ to $\Omega(t)$,
this leads us to reformulate the original problem on $\hat\Omega$ and introduce $\bPsi(t)$ (as a function of $\hat\Omega$) as an additional
solution field.
The evolution of $\bPsi(t)$ will be determined by the evolution of the boundary $\Gamma(t)$ of $\Omega(t)$ as
given by the free boundary problem and an harmonic extension into $\hat\Omega$.
In effect, we arrive at a formulation of the free boundary problem on the reference-domain using the domain transformation $\bPsi(t)$ as in Arbitrary-Lagrangian-Eulerian (ALE) methods \cite{hirt1974arbitrary,donea2004encyclopedia}.
The freedom in choice of the harmonic extension to determine $\bPsi(t)$ can be seen in
analogy to choosing a phase condition in context of the freezing formulation
discussed in \cite{OhlbergerRave2013}.

Literature on model order reduction for free boundary problems seems to be mostly non-existent.
We are only aware of the preliminary work in \cite{HinzeKrenciszekEtAl2014} and the following related works.
In the context of fluid-structure interaction (FSI) problems a similar ALE formulation has been considered
in \cite{BallarinRozza2016}.
The use of a reference domain for the reduction of models with parametrized geometry dates back to
the early days of reduced basis methods, e.g. \cite{PrudhommeRovasEtAl2001}.
This approach has been extended to FSI problems in \cite{LassilaQuarteroniEtAl2012}.
An increasingly popular way to deal with moving domains is based on an implicit description of the geometry through indicator functions (e.g. level sets \cite{sethian1999level,osher2006level}) as it often allows for a higher flexibility of the geometry handling w.r.t. large deformations and topology changes.
A nonlinear approximation method based on the truncation of functions in $V_N$ via time-dependent indicator functions
is discussed in \cite{BalajewiczFarhat2014}.
Due to the lack of hyperreduction no fully online-efficient reduced order model (ROM) is obtained, however.
Model order reduction of phase-field models, in which $\Gamma(t)$ is approximated by an easer-to-approximate diffuse interface layer,
is discussed in \cite{Volkwein2001,RedekerHaasdonk2015,GraessleHinze2017}.

Despite the high accuracy of projection-based ROMs, conservation properties of the original PDE model 
are usually only approximately preserved by the reduction process. 
The exact conservation of quantities such as the total mass of the system is often of particular interest, however.
In this work we will derive globally mass conserving ROMs based on a carefully chosen finite element
discretization of the free boundary problem in which the mass conservation constraint is implemented by testing the
variational formulation with a constant function.
Including the constant functions in the reduced space then ensures mass conservation of the Galerkin ROM.
To preserve this property under empirical interpolation \cite{barrault2004empirical}, we propose a rank-one
modification of the interpolated mass matrix to ensure exact yet efficient assembly of the constraint.

The inclusion of locally constant test functions in the reduced space to obtain locally mass conservative
flux reconstructions is considered in \cite{OhlbergerSchindler2015} to improve the efficiency of a localized a posteriori error
estimator in the context of localized model order reduction.
In \cite{CarlbergChoiEtAl2017} an alternative approach to preserve conservation properties is presented, which is based on the inclusion of
additional constraints in the least-squares Petrov-Galerkin minimization problem that is solved in the ROM.

\subsubsection*{Content and structure of the paper}

The paper is structured as follows.
We introduce a mathematical model for osmotic cell swelling as a model problem in Section~\ref{sec:osmo}. The Arbitrary-Lagrangian-Eulerian formulation of the full order model is then discussed in Section~\ref{sec:ALE} before we apply model order reduction in Section~\ref{sec:MOR}. Based on the numerical experiments in Section~\ref{sec:numex}, we discuss the performance and potential of this approach before we conclude.

\section{Mathematical model of an osmotic cell swelling problem}\label{sec:osmo}
As a model free boundary problem we consider a mathematical model of osmotic cell swelling that is also considered in a.o. \cite{lippoth2012classical,raetz2016diffuse,frischmuth1999numerical,zaal2012cell,zaal2013variational,zaal2015well}.
A membrane separates the interior of a cell which is filled with a fluid from the outside. Inside and outside of the cell a solute concentration is dissolved to which the membrane is impermeable. The outer concentration is assumed to be constant and known. An extension to the ``two-phase'' osmosis problem where the outer concentration field is also considered to be unknown can be found in \cite{lippoth2014stability,raetz2016diffuse}. We, however, restrict to the simpler ``one-phase'' case.
In this system, the membrane is subject to two acting forces: on the one hand, a surface tension force that only depends on the shape of the membrane and counteracts large curvatures; on the other hand, a force induced by the tendency to equilibrate the solute concentration across the membrane. The latter is modeled by Van't Hoff's law, which states that the pressure at the boundary is proportional to the concentration difference at the free boundary.

Let us denote by $u$ the solute concentration. Inside the cell, $u$ is subject to a linear unsteady diffusion equation with constant diffusion coefficient $\alpha$. Boundary conditions result naturally from the conservation of mass principle as the total solute concentration is constant.
Let $(0,T]$,~$T>0$ be the time interval of interest and $\Omega_0 \subset \mathbb{R}^d,~d=2,3$ be the initial domain of the cell with $\Gamma_0 := \partial \Omega_0$ the initial shape of the membrane. Then, the model for the solute concentration and the boundary motion is given by\begin{subequations} \label{eq:osmosisproblem}
\begin{align}
 \partial_t u - \alpha \Delta u & = 0 & \quad \text{ in } \Omega(t), \label{eq:diffusion}\\
 \qb u + \alpha \partial_{\n} u & = 0 & \quad \text{ on } \Gamma(t), \label{eq:masscons} \\
  - \beta \kappa + \gamma (u - \uext) & = \qb & \quad \text{ on } \Gamma(t).
                                                        \label{eq:vnmodel}
\end{align}
\end{subequations}
Here, $\qb$ is the velocity (in normal direction) of the boundary $\Gamma(t) := \partial \Omega(t)$. With the mean
curvature $\kappa$ (positive for convex domains) the first term $- \beta \kappa$ in \eqref{eq:vnmodel} models the effect of surface tension, whereas $\gamma (u-\uext)$ models the effect of the osmotic pressure. The constants $\beta$ and $\gamma$ are material constants depending only on the membrane and the solute. $\partial_{\n} u$ is the normal derivative of $u$ with $\n$ the outer normal to $\Omega(t)$. 

With Reynolds' transport theorem we easily see that \eqref{eq:masscons} implies conservation of the total solute:
\begin{equation} \label{eq:rey}
  \frac{d}{dt} \int_{\Omega(t)} u ~ d\bx = \int_{\Omega(t)} \underbrace{\partial_t u}_{=\alpha \Delta u}~ d\bx + \int_{\Gamma(t)} \qb u ~ d\bs
= \int_{\Gamma(t)} \alpha \partial_{\n} u + \qb u ~ d\bs = 0.
\end{equation}

\begin{figure}
    \centering
    \begin{tabular}{ccc}
    \includegraphics[width=0.3\textwidth, valign=c]{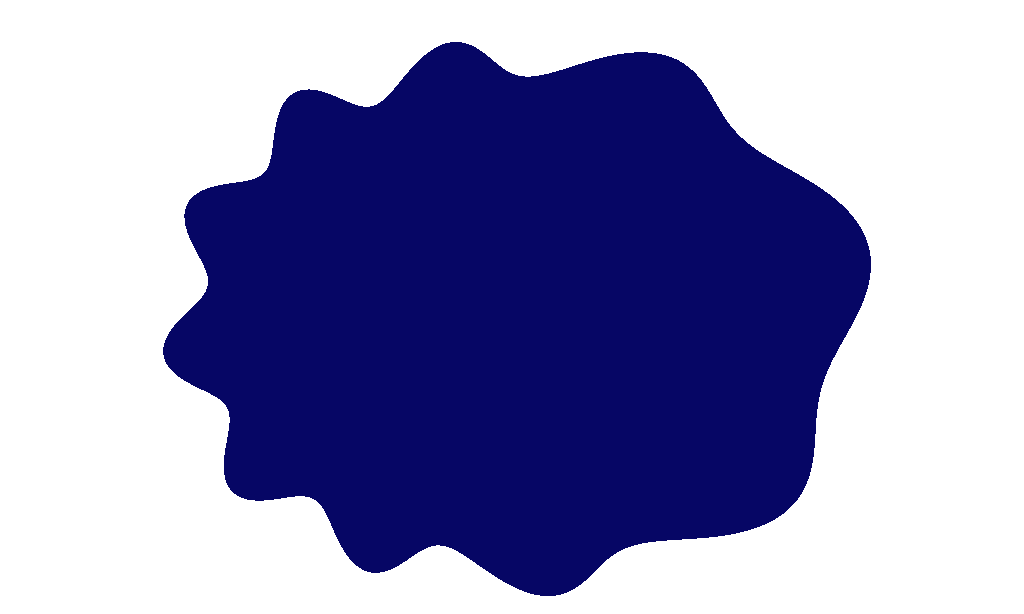}&
    \includegraphics[width=0.3\textwidth, valign=c]{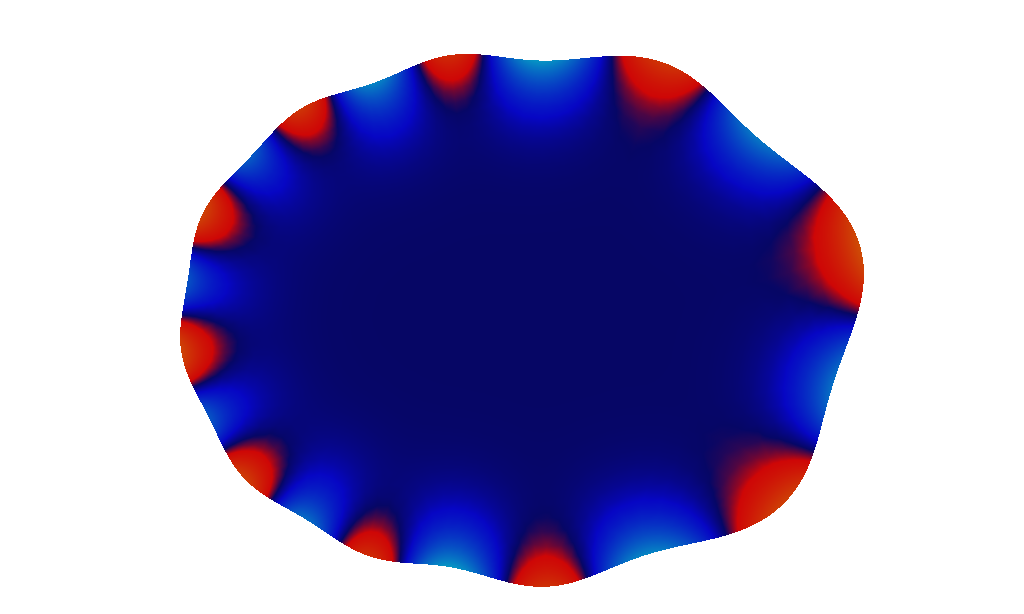}&
    \includegraphics[width=0.3\textwidth, valign=c]{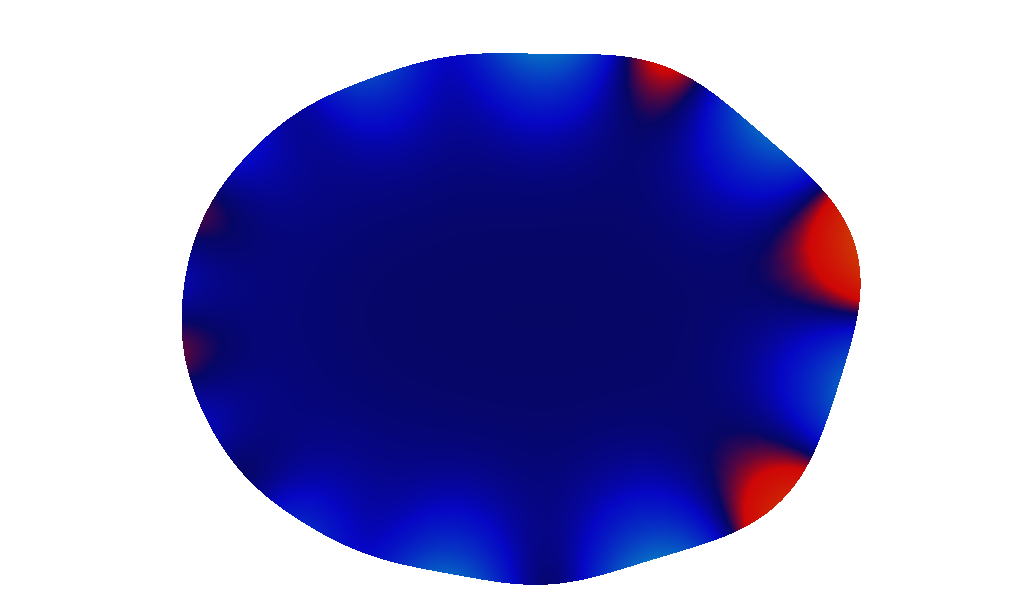}\\
    $t=0$ & $t = 0.25$ & $t = 0.5$ \\
    \includegraphics[width=0.3\textwidth, valign=c]{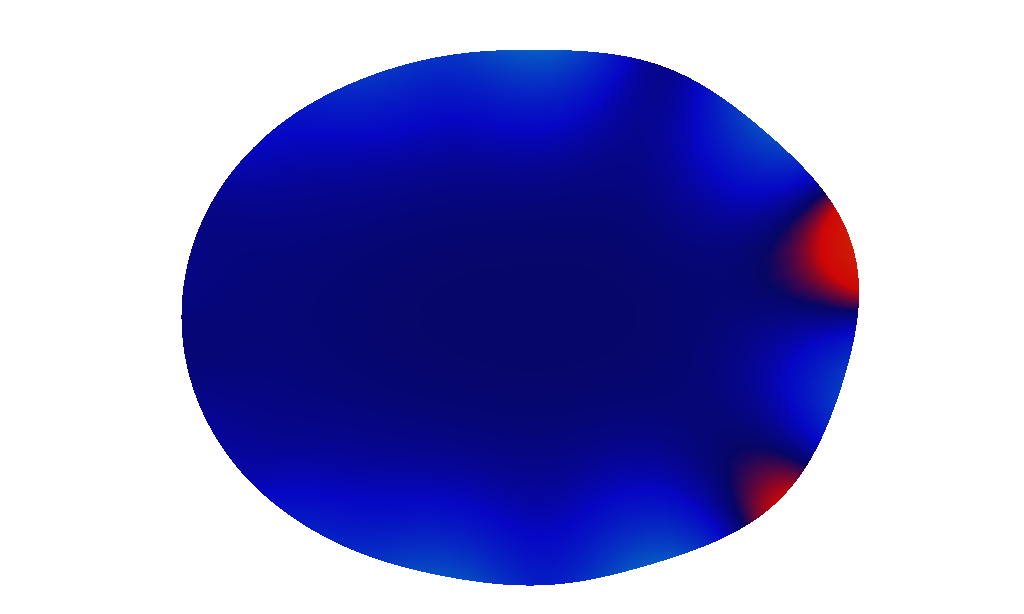}&
    \includegraphics[width=0.3\textwidth, valign=c]{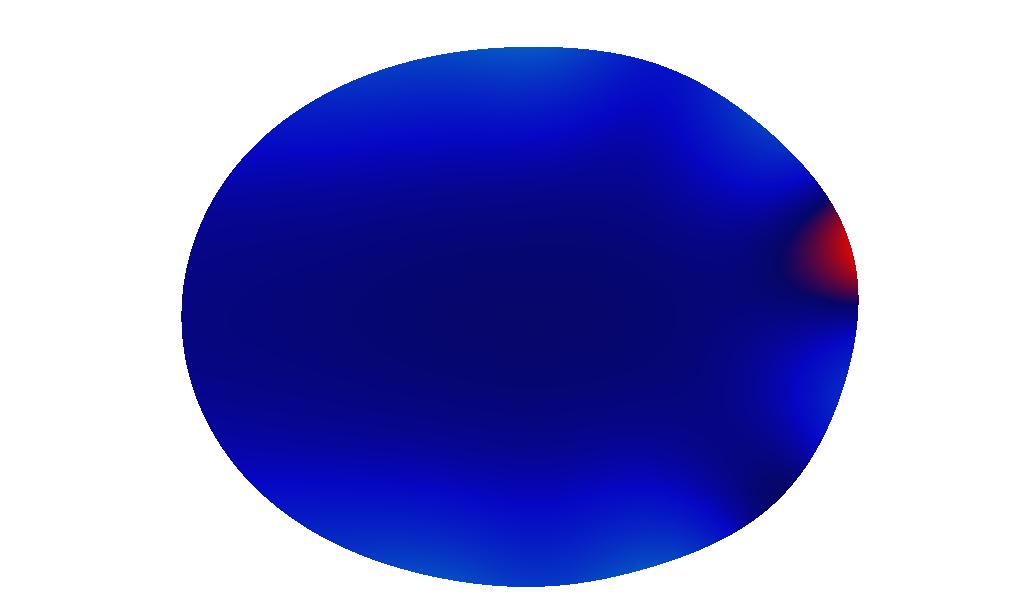}&
    \raisebox{-0.7cm}{\pgfplotscolorbardrawstandalone[
        colormap={osmosisc}{
            rgb(0.67374402284622203)=(0, 1, 1)
            rgb(0.95)=(0, 0, 1)
            rgb(1.0)=(0, 0, 0.50196078431400004)
            rgb(1.05)=(1, 0, 0)
            rgb(1.5269199609756501)=(1, 1, 0)
        },
        point meta min=0.67374402284622203,
        point meta max=1.5269199609756501,
        colorbar horizontal,
        colorbar style={xtick={0.67374402284622203,1.,1.5269199609756501}},
        parent axis width/.initial=0.25\textwidth
    ]} \\
    $t = 0.75$ & $ t = 1.0$ &
    \end{tabular}
    \caption{
      Evolution of cell geometry and concentration for an example considered in Section \ref{sec:numex:pmor}: The initially homogeneous distribution of the concentration is given in the initial cell shape. The mean curvature smoothens the cell shape, which leads to local changes in the concentration field. Finally, the evolution tends to a stationary state which is a circle with a homogeneous distribution of the concentration.
    }\label{fig:reconstruction}
\end{figure}

We notice that in the model a higher concentration at the boundary of the cell introduces a force that tends to expand the cell. The surface tension force typically has the opposite tendency. For convex domains it tends to compress the cell. Note that in this model the domain $\Omega$ cannot degenerate. This is due to the fact that the boundary condition \eqref{eq:masscons} ensures conservation of the total solute concentration so that the concentration increases if the cell shrinks. In the case of a shrinking cell, the concentration will eventually reach a level where the second term in \eqref{eq:vnmodel} compensates the surface tension force. 
For an initially simple connected domain $\Omega_0$ the system tends towards the stationary solution of \eqref{eq:osmosisproblem} which has a spherical shaped domain $\Omega_{\infty}$ with a constant concentration $u_\infty$. In Figure \ref{fig:reconstruction} an example evolution is shown.

\section{Arbitrary-Lagrangian-Eulerian-based finite element discretization} \label{sec:ALE}

In this section we derive a discretization method for the osmotic cell swelling problem based on an Arbitrary-Lagrangian-Eulerian (ALE) description \cite{hirt1974arbitrary,donea2004encyclopedia} and a finite element discretization. To this end, we first introduce the geometry description and a decomposition of the coupled problem into subproblems on the continuous level in Section \ref{sec:aleform}. In Section \ref{sec:prelimnote} we introduce some notation and a decoupling scheme consisting of discrete subproblems, the treatment of which is discussed one after another in Sections \ref{sec:laplacebeltrami} -- \ref{sec:concequation}. The conservation property
\eqref{eq:rey} is ensured on the discrete level in the discretization of the time stepping for the concentration field in Section \ref{sec:concequation}.

\subsection{ALE formulation of the continuous problem} \label{sec:aleform}
\subsubsection{Geometry description through mappings}
In an ALE description we deal with the fact that the domain is moving in time by introducing a
reference configuration $\hat\Omega$ and describe $\Omega(t)$ by a time-dependent transformation
\begin{equation}
  \bPsi: \hat\Omega \times [0,T] \to \mathbb{R}^d.
\end{equation}
By $\q = \partial_t \bPsi$
we denote the mesh velocity.
With abuse of notation, we write $\bPsi(t)$ for the restriction of $\bPsi$ to a fixed time $t$, i.e. the purely spatial function  $\bPsi(t): \hat{\Omega} \to \Omega(t)$.

\subsubsection{A generic extension from boundary transformation to volume transformation}
For the notation of restrictions of functions to the boundary $\Gamma$ we use a subindex, e.g. $\bPsi_\Gamma = \bPsi|_{\Gamma}$.
The evolution of the domain is determined by \eqref{eq:vnmodel},
\begin{equation}
  \partial_t \bPsi_{\Gamma} = \qb \ \n = -\beta \kappa \ \n + \gamma(u-\uext) \ \n \quad \text{ on } \Gamma(t).
\end{equation}
To determine the mesh
transformation $\bPsi(t)$ in the volume we use a linear extension operator $\hext : (H^{\frac12}(\Gamma(t)))^d \to (H^1(\Omega(t)))^d$ to extend the boundary transformation $\bPsi_{\Gamma}$ and the boundary velocity $\qb \ \n$ to the volume transformations $\bPsi$ and $\q$,
\begin{equation} \label{eq:velextension}
  \bPsi = \hext(\bPsi_{\Gamma})
    , \quad \Longrightarrow \quad \q = \partial_t \bPsi = \partial_t \hext(\bPsi_\Gamma) = \hext( \partial_t \bPsi_{\Gamma}) = \hext( \qb \ \n).
\end{equation}
\subsubsection{Continuous subproblems}
We write the ALE formulation for the continuous problem in terms of subproblems that are coupled
in a time interval $(t_{n-1},t_n]$. These subproblems are considered separately in the description of the discretization below:
\begin{subequations} \label{eq:coup}
  \begin{enumerate}
    \item
Fix a time $t=t_\ast$. Given $\bPsi(t_\ast)$ and $u(t_\ast)$, compute $\qb(t_\ast)$ as the solution to
  \begin{equation} \label{eq:coupa}
  \qb(t_\ast)  = - \beta \kappa + \gamma(u(t_\ast)-\uext) \quad  \text{ on }  \Gamma(t_\ast).
  \end{equation}
    \item
        Fix a time $t=t_\ast$. Given $\qb(t_\ast)$, compute a volumetric velocity as a suitable extension $\q(t_\ast)$:
  \begin{equation} \label{eq:coupb}
\q(t_\ast) = \hext(\qb(t_\ast) \n).
  \end{equation}
    \item
  Given the mesh velocity $\q(t)$ in a time interval $(t_{n-1},t_n]$ and the initial transformation at a time step $\Psiinit$, compute $\bPsi(t),~t\in (t_{n-1},t_n]$ as the solution to
    \begin{equation}  \label{eq:coupc}
      \partial_t \bPsi (t) = \q,~~t \in (t_{n-1},t_n],\qquad\bPsi(t_{n-1}) = \Psiinit.
    \end{equation}
    \item
Given the domain $\Omega(t) = \bPsi(t) ( \hat\Omega)$ in a time interval $(t_{n-1},t_n]$ and initial concentration data for this time interval $\uinit$, compute $u(t)$ as the solution to:
  \begin{equation} \label{eq:coupd}
    \left\{
      \begin{array}{rlll}
      \partial_t u - \alpha \Delta u & = 0 & \text{ in } \Omega(t),&~t\in (t_{n-1},t_n],  \\
      \qb u + \alpha \partial_{\n} u & = 0 & \text{ on } \Gamma(t),&~t\in (t_{n-1},t_n], \\
        u(t_{n-1}) & = \uinit & \text{ in } \Omega(t_{n-1}).&
      \end{array}
    \right.
  \end{equation}
  \end{enumerate}
\end{subequations}

\subsection{Preliminaries, notation and decoupling scheme} \label{sec:prelimnote}
In the discretization in this study we consider a simplicial triangulation $\T$ defining the reference domain $\hat\Omega$. On $\T$ we use the standard finite element space
\begin{equation}
  \Uh := \{ v \in C^0(\hat\Omega) \mid v|_T \in \mathcal{P}^k(T),~\forall~T\in\T \}
\end{equation}
with $\mathcal{P}^k(T)$ the space of polynomials of degree at most $k$.
Further, we define the vector valued space $\bUh := (\Uh)^d$ and its trace space $\bUhG := \bUh|_{\hat\Gamma}$ where $\hat\Gamma := \partial \hat{\Omega}$.
For the concentration field we introduce the unknown $\huh{} \in \Uh$ corresponding to a mapped function $\uh{} = \huh{} \circ \Psih{-1}$.
Correspondingly, we have the unknown transformation field $\bPsi$ that is approximated by $\Psih{}$ in $\bUh$.
The mesh velocity $\q$ is approximated on the reference domain with $\refqh{} \in \bUh$ so that $\qh{} = \refqh{} \circ \Psih{-1}$.
The vector-valued boundary velocity $\qb \ \n$ is approximated with $\qbh{} = \refqbh{} \circ \Psih{-1}$ with $\refqbh{}$ in the trace space $\bUhG$.

For the discretization in time we consider an equidistant decomposition of $(0,T]$ into $N$ time intervals and define time steps $t_i := i \Delta t$ with $\Delta t = T/N$.
We denote the discrete solutions at a time step $t_i$ by $\uh{i}$ ($\refuh{i}$), $\qh{i}$ ($\refqh{i}$), $\qbh{i}$ ($\refqbh{i}$) and $\Psih{i}$. Further, we introduce the notation $\Omega_h^i$ for the mapped domain $\Omega_h(t_i) := \Psih{i}(\hat\Omega)$ and accordingly define $\Gamma_h^i := \Psih{i}(\hat\Gamma)$.

In the remainder of this study, we consider a weakly coupled first order time integration scheme for the full discretization of the ALE formulation \eqref{eq:coupa}--\eqref{eq:coupd}. One time step in the scheme consists of the successive application of the following steps:
\begin{enumerate}
\item Given $\Psih{n-1}$, $\uh{n-1}$ compute $\refqbh{n-1}$ approximating \eqref{eq:coupa}, cf. Section  \ref{sec:laplacebeltrami}.
\item Extend $\refqbh{n-1}$ to $\hat\Omega$ resulting in $\refqh{n-1}$ , cf. Section \ref{sec:extension}.
\item With $\Psih{n-1}$ and $\refqh{n-1}$ given, we approximate \eqref{eq:coupc} with an explicit Euler step:
  \begin{equation}
    \Psih{n} = \Psih{n-1} + \Delta t \refqh{n-1}.
  \end{equation}
\item Take $\Psih{n-1}$ and $\Psih{n}$ to approximate the domain evolution. With $\uh{n-1}$ compute an approximation $\uh{n}$ to \eqref{eq:coupd}, cf. Section \ref{sec:concequation}.
\end{enumerate}
Below, for the application of integral transformations corresponding to a mapping $\bPsi$, we make use of the following notations for the Jacobian $\F$, the Jacobian determinant $J$ which is also the ratio of volume measure between preimage and image of $\bPsi$, the normal to the mapped domain $\n$, the ratio of the surface measures $J_\Gamma$ and the tangential projection onto the mapped domain $\P$:\vspace*{-0.15cm}
 \\
\begin{subequations} \label{eq:trafonotation}
\begin{minipage}{\textwidth}
\begin{minipage}{0.45\textwidth}
\begin{align}
    \F(\bPsi) &:= D\bPsi, \\
    J(\bPsi) &:= \det(\F(\bPsi)), \\
    J_\Gamma(\bPsi) &:= \|\F^{-T}(\bPsi)\hat\n\|\ J(\bPsi),
\end{align}
\end{minipage}
\hfill
\begin{minipage}{0.52\textwidth}
\begin{align}
    \n(\bPsi) &:= \F^{-T}(\bPsi)\hat\n \ \|\F^{-T}(\bPsi)\hat\n\|^{-1},\\
  \P(\bPsi) &:= \I - \n(\bPsi) \otimes \n^T(\bPsi). \\
   & \nonumber
\end{align}
\end{minipage}\vspace*{0.25cm}
\end{minipage}  
\end{subequations}
We notice that $\det(\F(\bPsi))$ is positive as long as $\bPsi$ is sufficiently close to the identity, i.e. as long as the deformation of the domain does not get too large.
We define $\F^i$ as the Jacobian to $\Psih{i}$, i.e. $\F^i = \F(\Psih{i})$ and define $J^i$, $J_\Gamma^i$, $\n^i$ and $\P^i$ accordingly.

\subsection{Discretization of the boundary velocity} \label{sec:laplacebeltrami}
We want to approximate $\qb \ \n$ with $\qb$ as in \eqref{eq:coupa} for $t_*=t_{n-1}$.
Hence, we seek for an approximation $\refqbh{n-1} \in \bUhG$ to
  \begin{equation} \label{eq:coupa2}
\refqbh{n-1} \approx \qb \ \n = - \beta \kappa \ \n + \gamma(u(t_\ast)-\uext) \ \n \quad  \text{ on }  \Gamma_h^{n-1}.
  \end{equation}
Multiplying with a test function $\bs_h = \hat{\bs}_h \circ \left(\Psih{n-1}\right)^{-1}$, $\hat{\bs}_h \in \bUhG$ and integrating over $\Gamma_h^{n-1}$ yields

\begin{align}
\int_{\Gamma_h^{n-1}} \qbh{n-1} \cdot \bs_h~d\bs 
  & =
- \beta \kappa_h^n (\bs_h) + \gamma \int_{\Gamma_h^{n-1}} (\uh{}-\uext)(\bs_h \cdot \n) ~ d\bs,
\end{align}
where $\kappa_h^n (\bs_h)$ is a discrete curvature linear form. 
To compute the mean curvature we make use of two main ideas from \cite{baensch}. First, we use the Laplace-Beltrami characterization of the mean curvature, $ - \kappa \ \n = \Delta_\Gamma \id $ in a weak formulation, to avoid the computation of second derivatives. This allows to make sense of a curvature even for polygonal boundaries.
Secondly, it is well-known that an explicit treatment of the curvature in free boundary problems leads to (severe) time step restrictions. These can be circumvented using an implicit approximation of the curvature. Hence, we aim at computing the curvature at time $t_{n}$ instead of $t_{n-1}$.
As $\Gamma_h^n$ is not known, we approximate the identity operator on the boundary to time $t_{n}$ by $\bx + \Delta t ~ \qbh{n-1} = \id|_{\Gamma_h^{n-1}} + \Delta t ~ \qbh{n-1} $ which yields the discrete curvature linear form
  \begin{equation}
\kappa_h^n (\bs_h) :=
 \int_{\Gamma_h^{n-1}} \nabla_{\Gamma} \id : \nabla_{\Gamma} \bs_h  ~ d\bs
 + \Delta t \int_{\Gamma_h^{n-1}} \nabla_{\Gamma} \qbh{n-1} : \nabla_{\Gamma} \bs_h ~ d\bs.
\end{equation}
We notice that the superscript $n$ at the linear form indicates that the curvature computation corresponds to time $t_n$.
Hence, our discretization of \eqref{eq:coupa} is: Find $\qbh{n-1} = \refqbh{n-1} \circ (\Psih{n-1})^{-1}$ with $\refqbh{n-1} \in \bUhG$, s.t.\ for all $\bs_h = \hat{\bs}_h \circ (\Psih{n-1})^{-1}$ with $\hat{\bs}_h \in \bUhG$ there holds
\begin{align}
\int_{\Gamma_h^{n-1}} \qbh{n-1} \cdot \bs_h~d\bs +& 
  \beta \Delta t \int_{\Gamma_h^{n-1}}  \nabla_{\Gamma} \qbh{n-1} : \nabla_{\Gamma} \bs_h~d\bs \\
  & =
- \beta \int_{\Gamma_h^{n-1}} \nabla_{\Gamma} \id : \nabla_{\Gamma} \bs_h ~ d\bs + \gamma \int_{\Gamma_h^{n-1}} (\uh{}-\uext)(\bs_h \cdot \n) ~ d\bs. \nonumber
\end{align}
Here, the tangential gradient acts on $\Gamma_h^{n-1}$, i.e. $\nabla_\Gamma \id = \P \cdot \nabla \id = \P$ where $\P = \I - \n \otimes \n^T$ is the tangential projection.
With the notations from \eqref{eq:trafonotation}  
and $\hat{\n}$ the outer normal to $\hat\Omega$, we can write this as an equation on $\hat{\Gamma}$:
Find $\refqbh{n-1} \in \bUhG$, s.t.\ for all $\hat{\bs}_h \in \bUhG$  there holds
\begin{equation} \label{eq:boundvel:ref}
  \begin{split}
  & \!\!\!\! \int_{\hat \Gamma} J_\Gamma^{n-1} \refqbh{n-1} \cdot \hat{\bs}_h \ d\bs
 + \beta \Delta t \int_{\hat \Gamma} J_\Gamma^{n-1}  \big( \P \cdot (\F^{n-1})^{-T} \cdot \nabla \refqbh{n-1} \big)  : \big( (\F^{n-1})^{-T} \nabla \hat{\bs}_h \big)~d\bs  \!\!\!\!\!\!\!\!\\
 \!\!\!\! = & - \beta \int_{\hat \Gamma} J_\Gamma^{n-1}  \P : (\F^{n-1})^{-T} \nabla_{\hat \Gamma} \ \hat{\bs}_h ~ d\bs + \gamma \int_{\hat \Gamma} J_\Gamma^{n-1}  (\huh{}-\uext)  \hat{\bs}_h \cdot ( (\F^{n-1})^{-T} \hn  ) \ d\bs. \!\!\!\!\!\!\!\!\!\!\!\!
\end{split}
\end{equation}

\subsection{Extension of the boundary velocity} \label{sec:extension}
For a given boundary velocity $\refqbh{n-1}$ (respectively $\qbh{n-1}$) we seek for an extension $\refqh{n-1} = \hexth(\refqbh{n-1}) \in \bUh$ (respectively $\qh{n-1}$) and use the solution operator (in a standard FEM formulation) of an harmonic extension (on the reference domain), i.e. we define $\refqh{n-1} = \hexth(\refqbh{n-1}) \in \bUh$ as the solution to
\begin{subequations} \label{eq:linelas}
  \begin{align}
      - \operatorname{div}[h_\mathcal{T}^{-1} ( \nabla \refqh{n-1} + (\nabla \refqh{n-1})^T)] & = 0 & \quad \text{ in } \hat\Omega, \\
    \refqh{n-1} & = \refqbh{n-1} & \quad \text{ on } \partial \hat\Omega,
  \end{align}
\end{subequations}
where $h_{\mathcal{T}}(\bx)$ is the locally constant grid function assigning to each $\bx \in T \in \mathcal{T}$ the diameter of $T$.

Note that one easily finds more sophisticated choices of the extension operator in the literature which allow to provide more control on the shape regularity of deformed mesh for larger deformations, see e.g. \cite[Section 5.1.2 and 5.1.3]{donea2004encyclopedia} and the references therein. In this study we consider only moderate deformations and take the liberty to consider only the simplified choice of an harmonic extension. We notice however that this restriction is not crucial for the applicability of the model order reduction considered below but simplifies the presentation as this extension operator is linear and parameter independent.

\subsection{Conservative concentration update on a moving domain} \label{sec:concequation}
We derive a time stepping procedure tailored to preserve the global solute mass. First, let us assume that a continuous mapping $\bPsi : \hat{\Omega} \times (t_{n-1},t_n] \to \mathbb{R}^d$ is known.
To $\hat{v} \in \bUh$, we define the mapped function $v(\bx,t) = \hat{v}(\bPsi^{-1}(t)(\bx))$ and apply Reynolds' transport theorem to the product $u\, v$ where $u$ is the exact solution to \eqref{eq:coupd}. This gives
\begin{align}
  \frac{d}{dt} \int_{\Omega(t)} u v ~ d\bx
  &= \int_{\Omega(t)} \partial_t(uv)~ d\bx + \int_{\partial \Omega} \qb u v ~ d\bs \nonumber \\
  &= \int_{\Omega(t)} \alpha \Delta u v - (\q \cdot \nabla v) u ~ d\bx + \int_{\partial \Omega} \qb u v ~ d\bs  \\
  &= \int_{\Omega(t)} - \alpha \nabla u \cdot \nabla v - (\q \cdot \nabla v) u ~ d\bx, \qquad t\in (t_{n-1},t_n]. \nonumber 
\end{align}
where we made use of $\partial_t u = \alpha \Delta u$ (from \eqref{eq:coupd}), $\partial_t v = - \q \cdot \nabla v$ (chain rule) and the boundary conditions in \eqref{eq:coupd}. Integration over $(t_{n-1},t_n]$ then yields
\begin{equation}
\int_{\Omega(t_n)} u v ~ d\bx - \int_{\Omega(t_{n-1})} u v ~ d\bx
  = \int_{t_{n-1}}^{t_n} \int_{\Omega(t)} - \alpha \nabla u \cdot \nabla v - (\q \cdot \nabla v) u ~ d\bx ~ dt.
\end{equation}
By choosing $v = 1$ ($\hat{v} = 1$) we recover the conservation of the total solute concentration  \eqref{eq:rey}. 
To arrive at a discretization, we replace the time integral with the right hand side rule, the exact geometries with $\Omega_h^{n-1}$ and $\Omega_h^{n}$ and the solution $u$ with the finite element approximations $\uh{n-1}$ and $\uh{n}$, respectively, yielding the discrete problem:

Find $\uh{n} = \huh{n} \circ \left(\Psih{n}\right)^{-1}$ with $\huh{n} \in \Uh$ s.t.\ for all $\vh = \hvh \circ \left(\Psih{n}\right)^{-1}$ with $\hvh \in \Uh$ there holds
\begin{equation} \label{eq:spacetime}
  \int_{\Omega_h^n} \uh{n} \vh~ d\bx + \Delta t \int_{\Omega_h^n} \uh{n}~\qh{n-1} \cdot \nabla \vh + \alpha \nabla \uh{n} \cdot \nabla \vh~ d\bx  = \int_{\Omega_h^{n-1}} \uh{n-1} \vh ~ d\bx.
\end{equation}
Equivalently, we can formulate the discretization with respect to the reference domain $\hat{\Omega}$:
Find $\huh{n} \in \Uh$, s.t.\ for all $\hvh \in \Uh$ there holds
\begin{align}\label{eq:conc:ref}
  \int_{\hat\Omega} J^{n}~\huh{n} \hvh \ d\bx
    & + \Delta t \int_{\hat \Omega} J^{n}~\huh{n}~\refqh{n-1} \cdot ((\F^{n})^{-T} \cdot \nabla \hvh)  d\bx \\
  & + \alpha~ \Delta t \int_{\hat \Omega} J^{n}~ ((\F^{n})^{-T} \nabla \huh{n}) \cdot ((\F^{n})^{-T} \nabla \hvh) \ d\bx
    = \int_{\hat \Omega} J^{n-1} ~ \huh{n-1} \hvh \ d\bx. \nonumber
\end{align}
We notice that the conservation property \eqref{eq:rey} is preserved also on the discrete level.

\section{Model order reduction} \label{sec:MOR}

\subsection{Parameterized full order model} \label{sec:param:model}

\begin{figure}
  \begin{center}
    \includegraphics[width=0.9\textwidth]{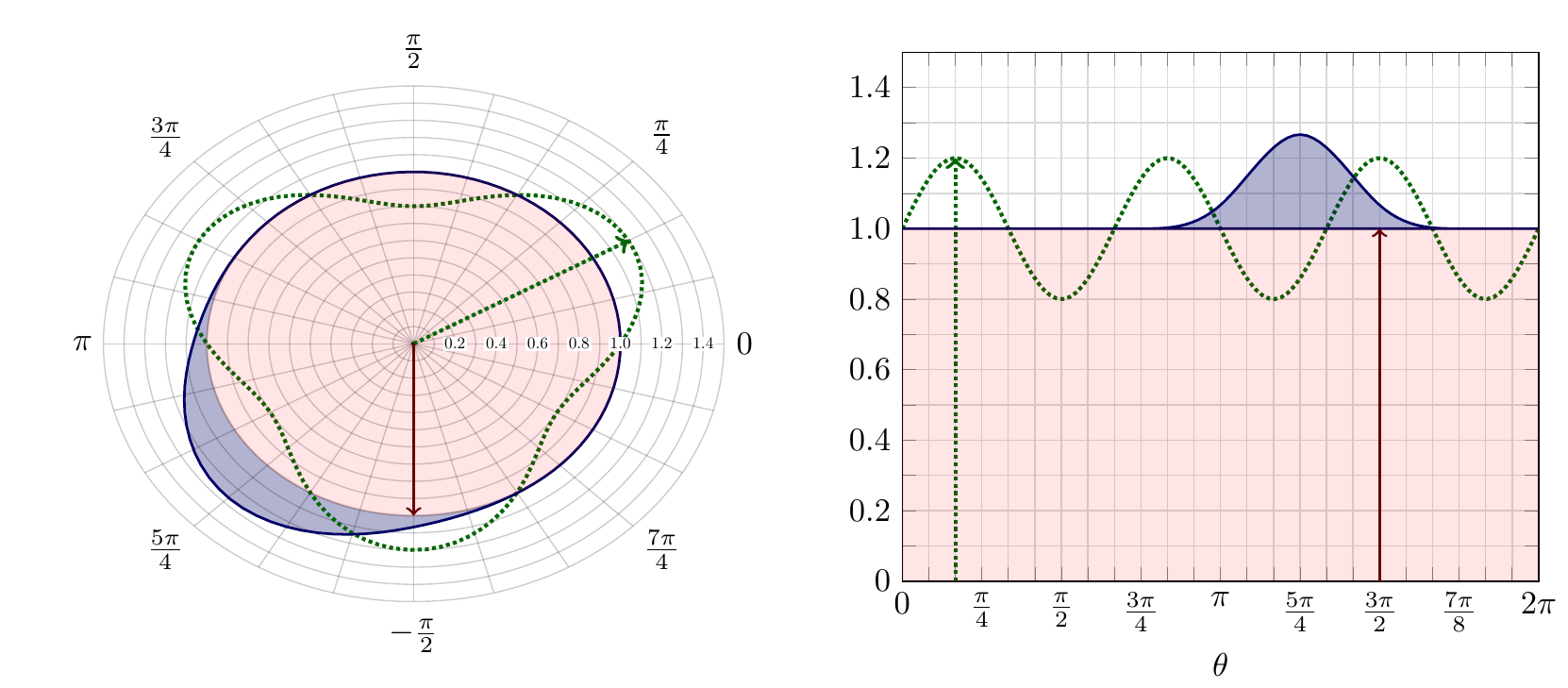} \vspace*{-0.35cm}
  \end{center}
  \caption{Sketch of possible initial shapes and corresponding graphs.}
  \label{fig:r}
\end{figure}

As parameters to the problem \eqref{eq:osmosisproblem} we consider the initial shape of the osmotic cell, the diffusivity $\alpha$ and the surface tension parameter $\beta$.
We notice that with dimensional analysis one easily checks that --- up to rescaling of time --- $\gamma$ is not an independent parameter in the model.
Thus, in the following we assume $\gamma$ to be constant.
Possible other parameters are the initial and the exterior concentration $\uinit=u(\cdot,0)$, $\uext$,
which we will not consider in this study, however.
For simplicity, we set $\uext=0$ in the sequel.
Regarding the parameterization of the initial domain $\Omega_0$,
we consider only initial shapes which are star-shaped and can be represented as a graph in normal direction, i.e.
\begin{equation}\label{eq:initialshape}
  \Gamma_0 = \{ \bx + r(\bx) \n(\bx) \mid \bx \in S_1(0) \},
\end{equation}
where $S_1$ is the unit sphere, $\n$ is the unit outer normal to $S_1$ and $r$ is $C^1$-smooth on $S_1$, cf.\ also Figure \ref{fig:r} for a sketch.
As reference domain we choose the unit disk $\hat\Omega:=D_2(0)$.
For simplicity we assume that $r$ is given as a linear combination of the form
\begin{equation}
    r(\bx) = \delta_1 r_1(\bx) + \ldots + \delta_L r_L(\bx).
\end{equation}
More general parameterizations could be considered using the empirical interpolation procedure described in
Section \ref{sec:ei}.
In total, the solution will depend on the parameter vector
\begin{equation}
    \mu := (\alpha, \beta, \delta_1, \ldots, \delta_L) \in \mathcal{P} \subset \mathbb{R}^{L+2},
\end{equation}
where $\mathcal{P}$ is the set of admissible parameters. 
The dependence on $\mu$ will be signified by adding the subscript $\mu$ to the respective solution fields.
The individual components of a parameter vector $\mu$ will be referred to as $\alpha_\mu, \beta_\mu, \delta_{1,\mu}, \ldots, \delta_{L,\mu}$.

Summarizing the discretization derived in Section~\ref{sec:ALE},
the discrete solution fields $\huhmu{n}$, $\Psihmu{n}$, $n = 0, \ldots N$ for $\mu \in \mathcal{P}$ are determined as follows:

As initial conditions we have
\begin{equation}\label{eq:fominitial}
    \huhmu{0} = \mathcal{I}_h(u_{init}), \qquad \Psihmu{0} = \id + \sum_{l=1}^L \delta_{l,\mu}\cdot\hexth(\mathcal{I}_{\Gamma,h}(r_l\,\n)),
\end{equation}
with $\mathcal{I}_h$, $\mathcal{I}_{\Gamma,h}$ denoting linear interpolation operators for $\bUh$, $\bUhG$
and $\hexth(\hbpsi_{\Gamma})$ denoting the solution of the extension problem defined in Section~\ref{sec:extension} for given
boundary data $\hbpsi_{\Gamma}$.

We introduce some notation for bilinear and linear forms appearing in
\eqref{eq:boundvel:ref} and \eqref{eq:conc:ref}
which simplify the presentation of the model order reduction approach below:
Problem \eqref{eq:boundvel:ref} becomes
\begin{subequations}\label{eq:fom}
\begin{gather}
\begin{multlined}[c][0.9\textwidth]
    \mqa(\refqbhmu{n-1}, \hat{\bs}_h; \Psihmu{n-1}) + \beta_\mu \Delta t \cdot \mqb(\refqbhmu{n-1}, \hat{\bs}_h; \Psihmu{n-1}) \\
    = - \beta_\mu \cdot \lqa(\hat{\bs}_h; \Psihmu{n-1}) +
    \gamma \cdot \lqc(\hat{\bs}_h; \Psihmu{n-1}, \huhmu{n-1}) \quad \forall \hat{\bs}_h \in  \bUhG.
\end{multlined}
\intertext{For the update of the transformation field we have}
\begin{multlined}
    \refqhmu{n-1} = \hexth(\refqbhmu{n-1}), \qquad \Psihmu{n} = \Psihmu{n-1} + \Delta t\, \refqhmu{n-1}
\end{multlined}
\intertext{where $\hexth$ is the extension operator from \eqref{eq:linelas}.
    The concentration update, cf. \eqref{eq:conc:ref}, reads as}
\begin{multlined}[c][0.9\textwidth]\label{eq:mordetcupdate}
    \muu(\huhmu{n}, \hvh; \Psihmu{n}) + \Delta t \cdot \aub(\huhmu{n}, \hvh; \Psihmu{n}, \refqhmu{n-1}) \\ + \alpha_\mu \Delta t \cdot \aua(\huhmu{n}, \hvh; \Psihmu{n})
    = \muu(\huhmu{n-1}, \hvh; \Psihmu{n-1}) \quad \forall \hat{v}_h \in  \Uh.
\end{multlined}
\end{gather}
\end{subequations}
In these equation systems, $\mqa, \mqb: H^{\frac{1}{2}}(\hat\Gamma)^d \times H^{\frac{1}{2}}(\hat\Gamma)^d \to \mathbb{R}$,
$\lqa, \lqc \in  (H^{\frac{1}{2}}(\hat\Gamma)^{d})^\prime$, $\muu, \aua, \aub: H^{\hat\Omega}\times H^{\hat\Omega} \to \mathbb{R}$,
are bilinear and linear forms given by

\begin{subequations}
\begin{align}
    \mqa(\hat{q}, \hat{\bs}; \hbpsi) &:= \int_{\hat{\Gamma}} \cmqa(\hbpsi) \hat{q} \cdot \hat{\bs} \ d\bs, &
    \aua(\hat{u}, \hat{v}; \hbpsi) &:= \int_{\hat\Omega} \caua(\hbpsi) \nabla\hat{u} \cdot \nabla\hat{v} \ d\bx, \\
    \mqb(\hat{q}, \hat{\bs}; \hbpsi) &:= \int_{\hat{\Gamma}} \cmqb(\hbpsi)\cdot \nabla_{\hat\Gamma} \hat{q} : \nabla_{\hat\Gamma}\hat{\bs} \ d\bs,\!\!\!\!\!\!\!& 
    \lqa(\hat{\bs}; \hbpsi) &:=  \int_{\hat\Gamma} \clqa(\hbpsi) : \nabla_{\hat\Gamma} \hat{\bs} \ d\bs, \\
    \muu(\hat{u}, \hat{v}; \hbpsi) &:= \int_{\hat\Omega} \cmuu(\hbpsi) \hat{u}~ \hat{v} \ d\bx, &
    \lqc(\hat{\bs}; \hbpsi, \hat{\varphi}) &:= \int_{\hat\Gamma} \clqc(\hbpsi, \hat\varphi) \cdot \hat{\bs} \ d\bs, \quad \\
    \aub(\hat{u}, \hat{v}; \hbpsi, \hbeta) &:= \int_{\hat\Omega} \caub(\hbpsi, \hbeta) \hat{u} \cdot \nabla\hat{v} \ d\bx, &
    &
\end{align}
\end{subequations}
with coefficient functions
\begin{subequations}\label{eq:fom_coeffs}
\begin{align}
    \cmqa(\hbpsi) &:= J_\Gamma(\hbpsi), &
    \caua(\hbpsi) &:= J(\hbpsi)\F^{-1}\!(\hbpsi)\F^{-T}\!(\hbpsi), \!\!\!\\
    \cmqb(\hbpsi) &:= J_\Gamma(\hbpsi) \F^{-1}\!(\hbpsi) \cdot \P(\hbpsi) \cdot \F^{-T}\!(\hbpsi), \!\!\!\!\!\!\!\!\!\! &
    \clqa(\hbpsi) &:=  J_\Gamma(\hbpsi) \F^{-1}\!(\hbpsi) \cdot \P(\hbpsi), \!\!\!\\
    \cmuu(\hbpsi) &:=  J(\hbpsi), &
    \clqc(\hbpsi, \hat\varphi) &:= J_\Gamma(\hbpsi) \F^{-T}\!(\hbpsi) \cdot \hat\n \cdot \hat\varphi, \!\!\!\\
    \caub(\hbpsi,\hbeta) &:= J(\hbpsi)\F^{-1}\!(\hbpsi)\cdot\hbeta,
\end{align}
\end{subequations}
where we made use of the notation from \eqref{eq:trafonotation}.

\subsection{Reduced basis approximation}

We construct a ROM for \eqref{eq:fominitial}, \eqref{eq:fom} via Galerkin projection onto
reduced order approximation spaces $\bUrG \subset \bUhG$, $\bUrPsi \subset \bUh$,
$\Uru \subset \Uh$ for the solution fields $\refqbhmu{n}$, $\Psihmu{n}$ and $\huhmu{n}$,~$n=0,..,N$.

Many different strategies have been discussed in the literature for constructing low-dimensional reduced approximation spaces
from solution snapshots of the full order model \eqref{eq:fominitial}, \eqref{eq:fom}.
In this study we choose a basic proper orthogonal decomposition (POD, \cite{Sirovich1987}) approach:

We assume that an appropriate finite set of training parameters $\mathcal{S}_{train} \subset \mathcal{P}$ has
been chosen and compute the snapshot sets
\begin{subequations}
\begin{align}
    \bM_{\Gamma}&:=\{\refqbhmu{n} \,|\, \mu \in \mathcal{S}_{train}, 0\leq n \leq N\}, \\
    \bM_{\Psi}&:=\{\Psihmu{n} - \id \,|\, \mu \in \mathcal{S}_{train}, 0\leq n \leq N\}, \\
    M_{u}&:=\{\huhmu{n} \,|\, \mu \in \mathcal{S}_{train}, 0\leq n \leq N\},
\end{align}
\end{subequations}
of the solution time trajectories of \eqref{eq:fominitial}, \eqref{eq:fom} for the parameter vectors $\mu \in \mathcal{S}_{train}$.
The affine shift in the definition of $\bM_\Psi$ is introduced in accordance with the
implementation used in Section~\ref{sec:numex}, where the mesh deformations $\Psihmu{k} - \id$
are used as solution variable.

The reduced approximation spaces $\bUrG,~\bUrPsi,~\Uru$ are now obtained from a POD
of $\bM_{\Gamma},~\bM_{\Psi},~M_{u}$, respectively, for a given relative truncation error tolerance $\varepsilon_{rb}$:
The reduced spaces are spanned by the first $K$ left-singular vectors of the linear mapping $\Phi_*$ sending the $k$-th
canonical basis vector of $\mathbb{R}^{|M_*|}$ to $s_k$ where
$\{s_1, \ldots, s_{|M_*|}\}$ is an arbitrary enumeration of the snapshot set.
The singular value decomposition is computed w.r.t.\ the $H^1(\hat\Omega)$-inner product on $\Uh$, the
$(H^1(\hat\Omega))^d$-inner product on $\bUh$, resp.\ the $(L^2(\hat\Gamma))^d$-inner product on $\bUhG$.
The truncation rank $K$ is determined as the minimal $K$ s.t. $\sigma_{K+1}(\Phi_*)/\sigma_1(\Phi_*) < \varepsilon_{rb}$, where
$\sigma_k(\Phi_*)$ denotes the $k$-th singular value of $\Phi_*$.

Denoting the $H^1(\hat\Omega)$-orthogonal projection of $\Uh$ onto $\Uru$ by $\Pru$ and the
$(H^1(\hat\Omega))^d$-orthogonal projection of $\bUh$ onto $\bUrPsi$ by $\bPrPsi$, the reduced
order model is given as follows:
For $\mu \in \mathcal{P}$, find $\hurmu{n} \in \Uru$, $\Psirmu{n} \in \id + \bUrPsi$, $0 \leq n \leq N$
with initial data
\begin{equation}\label{eq:rominitial}
    \hurmu{0} = \Pru(\mathcal{I}_h(u_{init})), \qquad \Psirmu{0} = \id + \sum_{l=1}^L \delta_{l}\cdot \bPrPsi(\hexth(\mathcal{I}_{\Gamma,h}(r_l))),
\end{equation}
with the boundary velocity $\refqbrmu{n-1} \in \bUrG$ given by
\begin{subequations}\label{eq:rom}
\begin{gather}\label{eq:rom_bnd}
\begin{multlined}[c][0.9\textwidth]
    \mqa(\refqbrmu{n-1}, \hat{\bs}_r; \Psirmu{n-1}) + \beta_\mu \Delta t \cdot \mqb(\refqbrmu{n-1}, \hat{\bs}_r; \Psirmu{n-1}) \\
    = - \beta_\mu \cdot \lqa(\hat{\bs}_r; \Psirmu{n-1}) + \gamma \cdot
    \lqc(\hat{\bs}_r; \Psirmu{n-1}, \hurmu{n-1})  \quad \forall \hat{\bs}_r \in  \bUrG,
\end{multlined}
\intertext{with the transformation field update given by}
\begin{multlined}\label{eq:rom_def}
    \refqrmu{n-1} = \bPrPsi(\hexth(\refqbrmu{n-1})), \qquad \Psirmu{n} = \Psirmu{n-1} + \Delta t \refqrmu{n-1}
\end{multlined}
\intertext{and with the concentration field update given by}
\begin{multlined}[c][0.9\textwidth]\label{eq:rom_u}
    \muu(\hurmu{n}, \hvr; \Psirmu{n}) + \Delta t \cdot \aub(\hurmu{n}, \hvr; \Psirmu{n}, \refqrmu{n-1}) \\ + \alpha_\mu \Delta t \cdot \aua(\hurmu{n}, \hvr; \Psirmu{n})
    = \muu(\hurmu{n-1}, \hvr; \Psirmu{n-1}) \quad \forall \hvr \in  \Uru.
\end{multlined}
\end{gather}
\end{subequations}

In order to solve \eqref{eq:rominitial}, \eqref{eq:rom}, we choose orthonormal bases of
$\Uru$, $\bUrPsi$, $\bUrG$.
We assemble the matrix of the linear operator $\bPrPsi\circ\hexth$ with respect to these bases,
as well as the coefficient vectors
of $\Pru(\mathcal{I}_h(u_{init}))$, and $\bPrPsi(\hexth(\mathcal{I}_{\Gamma,h}(r_l)))$
for the initial condition \eqref{eq:rominitial}.
In each time step, we then assemble the matrices and vectors of all
bilinear forms $a_*$ and linear forms $l_*$ appearing in \eqref{eq:rom_bnd}, \eqref{eq:rom_u}.
After that, the effort related to the computation of the basis coefficients of $\hurmu{n}$, $\Psirmu{n}$ is
\begin{equation}\label{eq:rom_effort}
    \mathcal{O}((\dim \bUrG)^3 + \dim \bUrPsi \cdot \dim \bUrG + (\dim \Uru)^3),
\end{equation}
where the summands correspond to the solution of \eqref{eq:rom_bnd}, \eqref{eq:rom_def} and \eqref{eq:rom_u}.
However, due to the dependence of the bilinear / linear forms on $\Psirmu{n}$, $\refqrmu{n}$, $\hurmu{n}$,
the matrix assembly has to be carried out in each time step, requiring substantial computational effort
proportional to $\dim \Uh$.
In the following section we use empirical interpolation to overcome this issue.

\subsection{Online efficient simulation via empirical interpolation}\label{sec:ei}

In order to achieve a fast assembly of \eqref{eq:rom_bnd}, \eqref{eq:rom_u} we use
empirical interpolation \cite{barrault2004empirical}
to approximate the coefficient functions \eqref{eq:fom_coeffs} by
linear combinations
\begin{subequations}
\begin{align}
    c_i(\hbpsi) &\approx \sum_{m=1}^{M_i} \theta_i^m(\hbpsi) c_i^m, &
    \theta_i^m(\hbpsi) &= [c_i(\hbpsi)(\bx_i^m)]_{k_i^m}, & i&\neq4,7\\
    c_4(\hbpsi,\hbeta) &\approx \sum_{m=1}^{M_4} \theta_4^m(\hbpsi,\hbeta) c_4^m, &
    \theta_4^m(\hbpsi,\hbeta) &= [c_4(\hbpsi,\hbeta)(\bx_4^m)]_{k_4^m},\\
    c_7(\hbpsi,\hat\varphi) &\approx \sum_{m=1}^{M_7} \theta_7^m(\hbpsi,\hat\varphi) c_7^m, &
    \theta_7^m(\hbpsi,\hat\varphi) &= [c_7(\hbpsi,\hat\varphi)(\bx_7^m)]_{k_7^m},
\end{align}
\end{subequations}
where for $i=1,\ldots,7$ the functions $c^m_i$ no longer depend on $\hbpsi$, $\hbeta$ or $\hat{\varphi}$,
$\bx_i^m \in \hat\Omega$ (resp. $\bx_i^m \in \hat\Gamma$) are interpolation points for $c_i$ and
$k_i^m$ are vector indices resp.\ matrix indices selecting a scalar component of
the evaluation of $c_i$ at $\bx_i^m$.

Defining the linear and bilinear forms
\begin{subequations}
\begin{align}
    a_1^m(\hat{q}, \hat{\bs}) &:= \int_{\hat{\Gamma}} c_1^m \hat{q} \cdot \hat{\bs} \ d\bs, &
    a_5^m(\hat{u}, \hat{v}) &:= \int_{\hat\Omega} c_5^m \nabla\hat{u} \cdot \nabla\hat{v} \ d\bx, \\
    a_2^m(\hat{q}, \hat{\bs}) &:= \int_{\hat{\Gamma}} a_2^m\cdot \nabla_{\hat\Gamma} \hat{q} : \nabla_{\hat\Gamma}\hat{\bs} \ d\bs,&
    l_1^m(\hat{\bs}) &:=  \int_{\hat\Gamma} c_6^m : \nabla_{\hat\Gamma} \hat{\bs} \ d\bs, \\
    a_3^m(\hat{u}, \hat{v}) &:= \int_{\hat\Omega} c_3^m \hat{u} \hat{v} \ d\bx, &
    l_2^m(\hat{\bs}) &:= \int_{\hat\Gamma} c_7^m \cdot \hat{\bs} \ d\bs, \\
    a_4^m(\hat{u}, \hat{v}) &:= \int_{\hat\Omega} c_4^m \hat{u} \cdot \nabla\hat{v} \ d\bx, 
\end{align}
\end{subequations}
we obtain approximations $a_i \approx \sum_{m=1}^{M_i} \theta_i^m \cdot a_i^m$ and $l_i \approx \sum_{m=1}^{M_i} \theta_{i+5}^m
\cdot l_i^m$ yielding the update equations
\begin{subequations}\label{eq:romei}
\begin{align}\label{eq:rom_bnd_ei}
    \sum_{m=1}^{M_1}& \theta_1^m(\Psirmu{n-1})\cdot a_1^m(\refqbrmu{n-1}, \hat{\bs}_r) +
    \beta_\mu \Delta t \cdot \sum_{m=1}^{M_2} \theta_2^m(\Psirmu{n-1}) \cdot a_2^m(\refqbrmu{n-1}, \hat{\bs}_r) \\
    & = -\beta_\mu \cdot \sum_{m=1}^{M_6}\theta_6^m(\Psirmu{n-1}) \cdot l_1^m(\hat{\bs}_r) + \gamma \cdot
    \sum_{m=1}^{M_7} \theta_7^m(\Psirmu{n-1}, \hurmu{n-1})\cdot l_7^m(\hat{\bs}_r)  \quad \forall \hat{\bs}_r \in  \bUrG, \nonumber
\intertext{and}
    \label{eq:rom_u_ei}
    \sum_{m=1}^{M_3}& \theta_3^m(\Psirmu{n})\cdot a_3^m(\hurmu{n}, \hvr) +
    \Delta t \cdot \sum_{m=1}^{M_4} \theta_4^m(\Psirmu{n},\refqrmu{n-1}) \cdot a_4^m(\hurmu{n}, \hvr) \\ 
    & + \alpha_\mu \Delta t \cdot \sum_{m=1}^{M_5} \theta_5^m(\Psirmu{n}) \cdot a_5^m(\hurmu{n}, \hvr)
    = \sum_{m=1}^{M_3} \theta_3^m(\Psirmu{n-1}) \cdot a_3^m(\hurmu{n-1}, \hvr) \quad \forall \hvr \in  \Uru. \nonumber
\end{align}
\end{subequations}
After pre-assembly of the matrices of $a_i^m$ and coefficient vectors of $l_i^m$, the effort for the assembly
of the equation systems \eqref{eq:rom_bnd_ei}, \eqref{eq:rom_u_ei} for an arbitrary $\mu \in \mathcal{P}$ 
is of order
\begin{equation}\label{eq:rom_ei_effort}
    \mathcal{O}( M_{1,2} \cdot (\dim \bUrG)^2 +  M_{6,7} \cdot \dim \bUrG
    + M_{3,4,5}\cdot (\dim \Uru)^2),
  \end{equation}
  with $M_{1,2}:= M_1 + M_2$,~~$M_{3,4,5} := M_3 + M_4 + M_5$ and $M_{6,7} := M_6 + M_7$.

\paragraph{Computation of the interpolation data.}

We consider training sets of function evaluations
\begin{align}
    \mathcal{C}_i \!&:=\! \{\!\,[c_i(\Psihmu{n})(\bx_j)]_j \,|\, \mu \in \mathcal{S}_{train}, 0 \!\leq\! n \!\leq\! N, 1 \!\leq\! j \!\leq\! |\mathcal{X}_i| \} \!\subset\!
    \ell^\infty(\mathcal{X}_i)^{d_i},~i \neq 4,7, \nonumber \\
    \mathcal{C}_4 \!&:=\! \{\!\,[c_4(\Psihmu{n},\refqrmu{n-1})(\bx_j)]_j \,|\, \mu \in \mathcal{S}_{train}, 0 \!\leq\! n \!\leq\! N, 1 \!\leq\! j \!\leq\! |\mathcal{X}_i| \} \!\subset\!
    \ell^\infty(\mathcal{X}_i)^{d_4}, \!\!\!\!\!\!\!\!\!\!\!\!\\
    \mathcal{C}_7 \!&:=\! \{\!\,[c_7(\Psihmu{n},\huhmu{n})(\bx_j)]_j \,|\, \mu \in \mathcal{S}_{train}, 0 \!\leq\! n \!\leq\! N, 1 \!\leq\!j \!\leq\! |\mathcal{X}_i| \} \!\subset\!
    \ell^\infty(\mathcal{X}_i)^{d_7},  \nonumber 
\end{align}
where $\mathcal{X}_i = \{\bx_1, \ldots,\bx_{|\mathcal{X}_i|}\}$ is an appropriate finite subset of
$\hat\Gamma$ ($i = 1,2,6,7$) resp. $\hat\Omega$ ($i = 3,4,5$) and $d_i$ corresponds to the
shape of the values of $c_i$, i.e. $d_1=d_3=1$, $d_4=d_7=d$ and $d_2=d_5=d_6=d^2$.
Using the $\mathcal{C}_i$ as input for the greedy algorithm from \cite{barrault2004empirical},
interpreting tensor fields in $\ell^\infty(\mathcal{X}_i)^{d_i} \cong \ell^\infty(\dot\bigcup_{k=1}^{d_i}
\mathcal{X}_i)$ as scalar functions,
we obtain the desired interpolation basis functions $c_i^m \in \operatorname{span}(\mathcal{C}_i)$
and interpolation points $\bx_i^m \in \mathcal{X}_i$, $1 \leq k_i^m \leq d_i$,
which approximate all elements of $\mathcal{C}_i$ with a relative $\ell^\infty(\mathcal{X}_i)^{d_i}$-error smaller
than a prescribed tolerance $\varepsilon_{ei}$.

Note that for the evaluation of the bilinear forms $a_i^m$ and the linear forms $l_i^m$, knowledge
of the corresponding coefficient function $c_i^m$ is only required at the finitely many quadrature
points used in the discretization scheme.
Hence, it is sufficient to choose $\mathcal{X}_i$ as the set of all these quadrature points.
Also note that our approach differs from \cite{BallarinRozza2016} in which each coefficient tensor component
is interpolated separately, whereas we perform a single interpolation of the full tensor field using
tensor-valued interpolation basis functions and scalar components of the tensor at given $\bx \in \mathcal{X}_i$
as interpolation points.
Since separate empirical interpolation of the tensor field components will in general select different interpolation
points $\bx \in \mathcal{X}_i$ for the individual components, we expect our approach to be more efficient in general
(i.e.\ require less coefficient tensor evaluations).

\paragraph{Minimally intrusive implementation of empirical interpolation.}
In many cases it is technically difficult to replace in the PDE solver's matrix assembly code the analytically defined coefficient functions $c_i$ by a vector of function evaluations
$c_i^m$ at the given quadrature points.
In the numerical example in Section \ref{sec:numex} we have used the following less intrusive approach to
implement the empirical interpolation procedure, which in addition 
does not require knowledge of the exact quadrature points used by the assembly routine:

Noting that $c_i^m \in \operatorname{span}(\mathcal{C}_i)$, we can represent $c_i^m$, $i = 1,2,3,5$, as
   \begin{equation}
       c_i^m = \sum_{\mu \in \mathcal{S}_{train}}\,\sum_{n=1}^N \gamma_{i,\mu}^{m,n}\cdot [c_i(\Psihmu{n})(\bx_j)]_j,
   \end{equation}
where the linear coefficients $\gamma_{i,\mu}^{m,n}$ can be directly obtained from the execution of the greedy algorithm.
Using this representation, it immediately follows that
   \begin{equation}
       a_i^m(\cdot, \cdot) =
       \sum_{\mu \in \mathcal{S}_{train}}\,\sum_{n=1}^N \gamma_{i,\mu}^{m,n}\cdot a_i(\cdot, \cdot; \Psihmu{n}).
   \end{equation}
Since the matrices $a_i(\cdot, \cdot; \Psihmu{n})$ have already been computed by the PDE solver,
we can easily assemble the matrix of $a_i^m$ using this formula.
The same argument applies to $\aub$, $l_1^m$, $l_2^m$.

\subsection{Global mass conservation}

As discussed in Section \ref{sec:concequation}, choosing the test function $\hvh \equiv 1$ in the concentration update equation
\eqref{eq:conc:ref} shows that the total mass $\int_{\Omega_h^n} \uh{n} d\bx$
is conserved by the discretization scheme, i.e. $\int_{\Omega_h^n} \uh{n} d\bx = \int_{\Omega_h^0} \uh{0} d\bx$ for all
$0 \leq n \leq N$.
By including the constant functions in the reduced concentration space $\Uru$,
\begin{equation}
    1 \in \Uru,
\end{equation}
the same argument can be applied to the reduced concentration field given by \eqref{eq:rom_u}:
\begin{align}
    \int_{\Omega_{r,\mu}^n} \urmu{n} d\bx &=
    \muu(\hurmu{n}, 1; \Psirmu{n}) + \Delta t \cdot \underbrace{\aub(\hurmu{n}, 1; \Psihmu{n}, \refqrmu{n-1})}_{= 0} + \nonumber
    \alpha_\mu \Delta t \cdot \underbrace{\aua(\hurmu{n}, 1; \Psirmu{n})}_{= 0} \\ &= \muu(\hurmu{n-1}, 1; \Psirmu{n-1})
     = \int_{\Omega_{r,\mu}^{n-1}} \urmu{n-1} d\bx. \label{eq:rom_u_cons}
\end{align}
Here, $\Omega_{r,\mu}^n := \Psirmu{n}(\hat\Omega)$ is the deformed domain at time step $n$ given by the reduced
transformation field $\Psirmu{n}$ and $\urmu{n} := \hurmu{n} \circ (\Psirmu{n})^{-1}$ the reduced concentration field on this
domain.

Note, however, that the argument in \eqref{eq:rom_u_cons} is no longer valid when replacing the bilinear forms $a_3$, $a_4$, $a_5$
by their empirical interpolants \eqref{eq:rom_u_ei}: While it still holds that
\begin{equation}
    \sum_{m=1}^{M_4} \theta_4^m(\Psirmu{n},\refqrmu{n-1}) \cdot \underbrace{a_4^m(\hurmu{n}, 1)}_{=0} =
    \sum_{m=1}^{M_5} \theta_5^m(\Psirmu{n}) \cdot \underbrace{a_5^m(\hurmu{n}, 1)}_{=0} = 0,
\end{equation}
the total mass at time step $n$ is only approximately given by the empirical interpolant of $a_3$:
\begin{equation}
    \int_{\Omega_{r,\mu}^n} \urmu{n} d\bx \approx
    \sum_{m=1}^{M_3} \theta_3^m(\Psirmu{n})\cdot a_3^m(\hurmu{n}, 1).
\end{equation}

One way to recover exact mass conservation in the ROM is to apply the empirical
interpolation procedure only to the coefficient functions of $a_4$, $a_5$, whereas 
$a_3(\cdot, \cdot; \Psirmu{n})$ is evaluated exactly.
For $d = 2$, we can evaluate $a_3$ as
\begin{equation}
    a_3(\cdot, \cdot, \Psirmu{n}) = \bar{a}_3(\cdot, \cdot, \Psirmu{n}, \Psirmu{n}),
\end{equation}
where $\bar{a}_3$ is the 4-tensor given by
\begin{equation}
    \bar{a}_3(\hat{u}, \hat{v}, \hbpsi, \hbeta) =
    \int_{\hat\Omega} (\partial_1 \hbpsi_1\cdot \partial_2 \hbeta_2 - \partial_1 \hbpsi_2\cdot \partial_2 \hbeta_1)
    \cdot \hat{u} \cdot \hat{v} \ d\bx.
\end{equation}
The effort to assemble the matrix of $a_3(\cdot, \cdot; \Psirmu{n})$ from the coefficients of $\bar{a}_3$
w.r.t.~a basis of $\Uru$ is of order $\mathcal{O}((\dim \Uru)^2\cdot(\dim \bUrPsi)^2)$.
Although highly efficient implementations for this operation are available, the higher computational
complexity in comparison to \eqref{eq:rom_effort}, \eqref{eq:rom_ei_effort} will lead to dominating runtime costs for large
reduced space dimensions.
In three spatial dimensions ($d=3$), the matrix of $a_3$ can be assembled exactly by the same argument
with a computational effort of $\mathcal{O}((\dim \Uru)^2\cdot(\dim \bUrPsi)^3)$.
We expect this to be non-favorable, even for relatively small dimensions of $\Uru$.
However, to ensure mass conservation, only the functional
$a_3(\cdot, 1; \Psirmu{n})$ needs to be known exactly, the matrix of which
can be computed by the same argument with a reduced effort of $\mathcal{O}((\dim \Uru)\cdot(\dim \bUrPsi)^d)$.

Thus, choosing a basis $\hat\varphi_i$, $1 \leq i \leq \dim \Uru$ for $\Uru$ such
that $\varphi_1 = 1$, we define the reduced bilinear form $\tilde{a}_3(\hat{u}_r, \hat{v}_r; \hbpsi_r)$
on $\Uru$ by
\begin{equation}
\tilde{a}_3(\hat{u}_r, \hat\varphi_i; \hbpsi_r):=
    \begin{cases}
        a_3(\hat{u}_r, 1; \hbpsi_r) & i = 1 \\
        \sum_{m=1}^{M_3} \theta_3^m(\hbpsi_r)\cdot a_3^m(\hat{u}_r, \hat\varphi_i) & 2 \leq i \leq \dim \Uru,
    \end{cases}
\end{equation}
and use as exactly mass conservative concentration update equation:
\begin{multline}\label{eq:rom_u_ei_cons}
    \tilde{a}_3(\hurmu{n}, \hvr; \Psirmu{n}) +
    \Delta t \cdot \sum_{m=1}^{M_4} \theta_4^m(\Psirmu{n},\refqrmu{n-1}) \cdot a_4^m(\hurmu{n}, \hvr) \\
    + \alpha_\mu \Delta t \cdot \sum_{m=1}^{M_5} \theta_5^m(\Psirmu{n}) \cdot a_5^m(\hurmu{n}, \hvr)
    = \tilde{a}_3(\hurmu{n-1}, \hvr; \Psirmu{n-1}) \quad \forall \hvr \in  \Uru.
\end{multline}
In total, this equation system can then be assembled with an effort of
\begin{equation}\label{eq:rom_ei_effortb}
    \mathcal{O}( M_{3,4,5} \cdot (\dim \Uru)^2 + (\dim \Uru)\cdot(\dim \bUrPsi)^d).
\end{equation}

\section{Numerical experiments} \label{sec:numex}

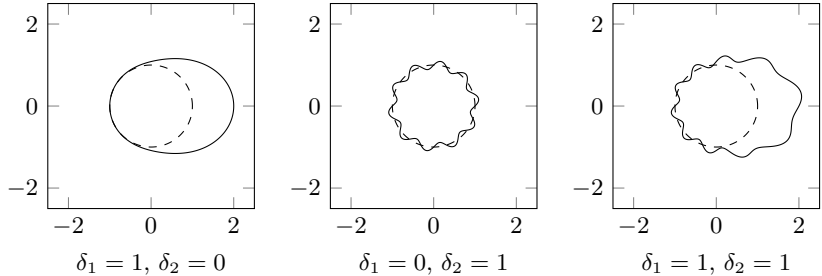
\begin{figure}[t] \centering
 \begin{tikzpicture}
     \begin{groupplot}[
             group style={group size=3 by 1,
                          horizontal sep=1.0cm,
                          vertical sep=1.5cm,},
             width=4.3cm,
             height=4.3cm,
             /tikz/font=\small,
             xmin=-2.5,
             xmax=2.5,
             ymin=-2.5,
             ymax=2.5,
     ]
         \nextgroupplot[xlabel={$\delta_1 = 1$, $\delta_2 = 0$}]
            \addplot[samples=360, domain=0:360, dashed] ({cos(x)}, {sin(x)});
            \addplot[samples=360, domain=-180:180] ({cos(x) * (1 + exp(- (x/180*pi)^2))},
                                                    {sin(x) * (1 + exp(- (x/180*pi)^2))});
         \nextgroupplot[xlabel={$\delta_1 = 0$, $\delta_2 = 1$}]
            \addplot[samples=360, domain=0:360, dashed] ({cos(x)}, {sin(x)});
            \addplot[samples=360, domain=-180:180] ({cos(x) * (1 + 0.1 * sin(10*x))},
                                                    {sin(x) * (1 + 0.1 * sin(10*x))});
         \nextgroupplot[xlabel={$\delta_1 = 1$, $\delta_2 = 1$}]
            \addplot[samples=360, domain=0:360, dashed] ({cos(x)}, {sin(x)});
            \addplot[samples=360, domain=-180:180] ({cos(x) * (1 + exp(- (x/180*pi)^2) + 0.1 * sin(10*x))},
                                                    {sin(x) * (1 + exp(- (x/180*pi)^2) + 0.1 * sin(10*x))});
     \end{groupplot}
\end{tikzpicture}
\caption{Initial boundary $\Gamma_0$ for different parameter combinations $\delta_1$, $\delta_2$ in the
numerical experiment (solid line).}\label{fig:initial}
\end{figure}

As a test for the developed reduced order modeling workflow we consider the parameterized model from Section
\ref{sec:param:model} in two spatial dimensions and a two-dimensional parameterization ($L = 2$) of the initial boundary
$\Gamma_0$ given by
\begin{subequations}
\begin{align}
    r_1(\bx) &:= e^{-\theta(\bx)^2} \cdot \bx, &\theta(\bx) &:= \operatorname{atan2}(\bx),\\
    r_2(\bx) &:= 10^{-1} \sin(10\cdot\theta(\bx)) \cdot \bx,
\end{align}    
\end{subequations}
where $\operatorname{atan2}(\bx)$ denotes the $(-\pi, \pi]$-valued angle between $\bx$ and the
$
\begin{bmatrix}
    1 & 0 
\end{bmatrix}
$-axis (cf.\ Figure \ref{fig:initial}).
We choose $\gamma = 0.1$ and $T = 1$ as final simulation time.
The resulting ALE formulation is discretized using a first-order finite element approximation with
3,988 degrees of freedom per spatial variable, i.e. 11,964 degrees of freedom in total.
As time step size we choose $\Delta t = 0.01$.

\begin{figure}[p]
\centering
\begin{tabular}{cccc}
    $t$ & $\huh{n}$ & $(\Psih{n} - \id)_1$ & $(\Psih{n} - \id)_2$ \\
    $0.00$ &
    \includegraphics[width=0.25\textwidth, valign=c]{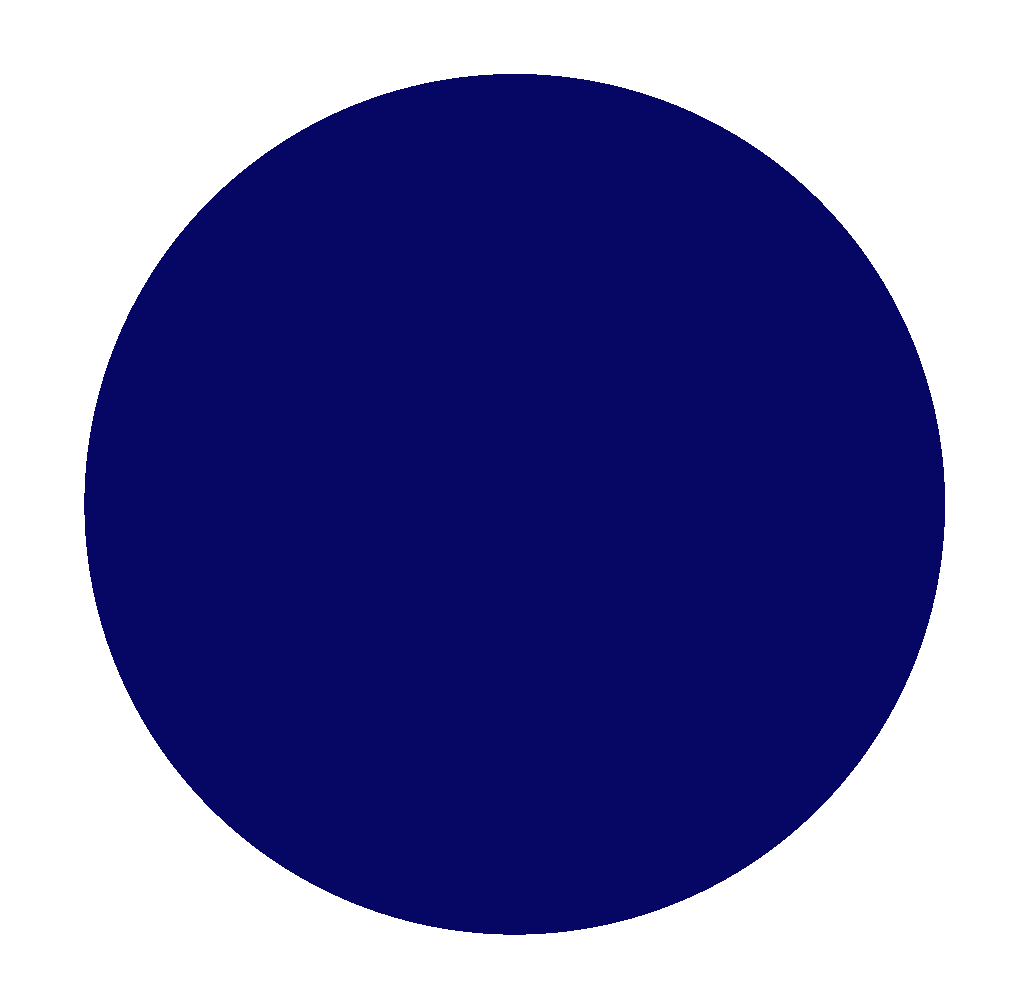} &
    \includegraphics[width=0.25\textwidth, valign=c]{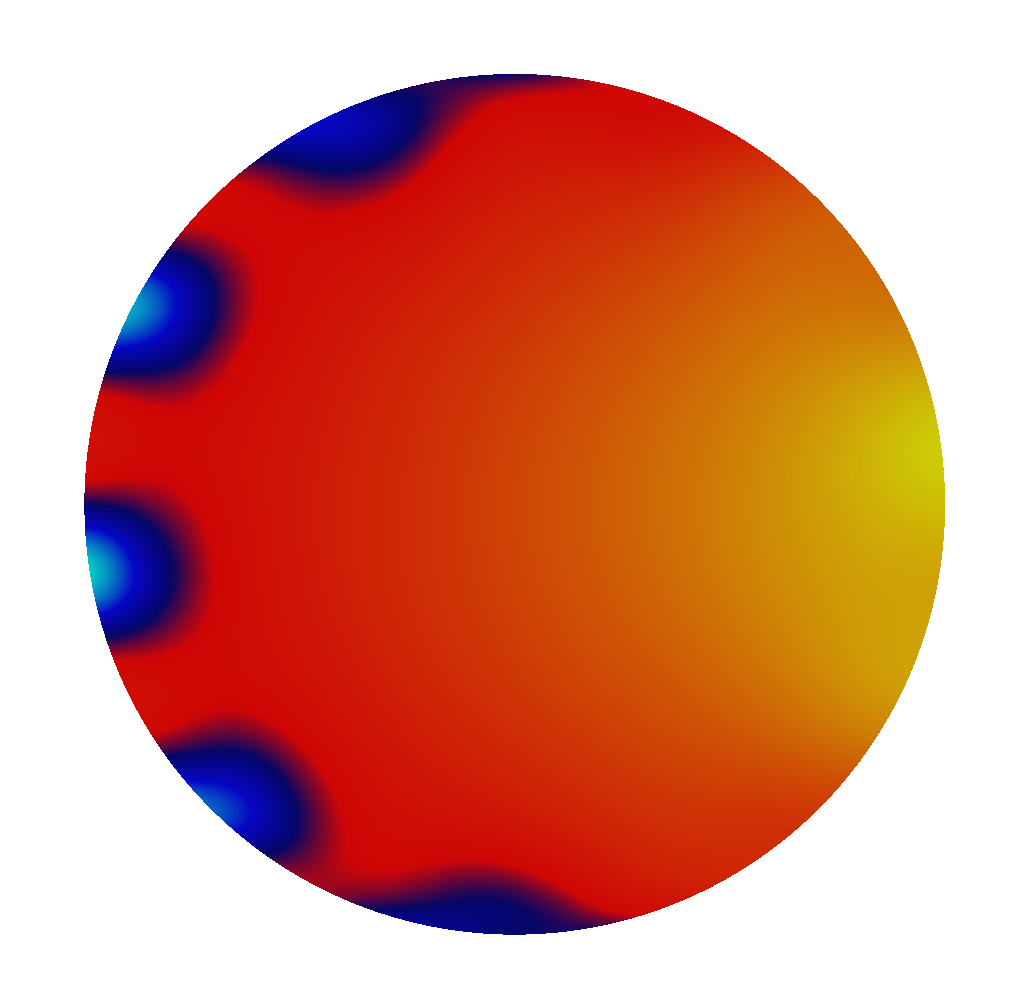} &
    \includegraphics[width=0.25\textwidth, valign=c]{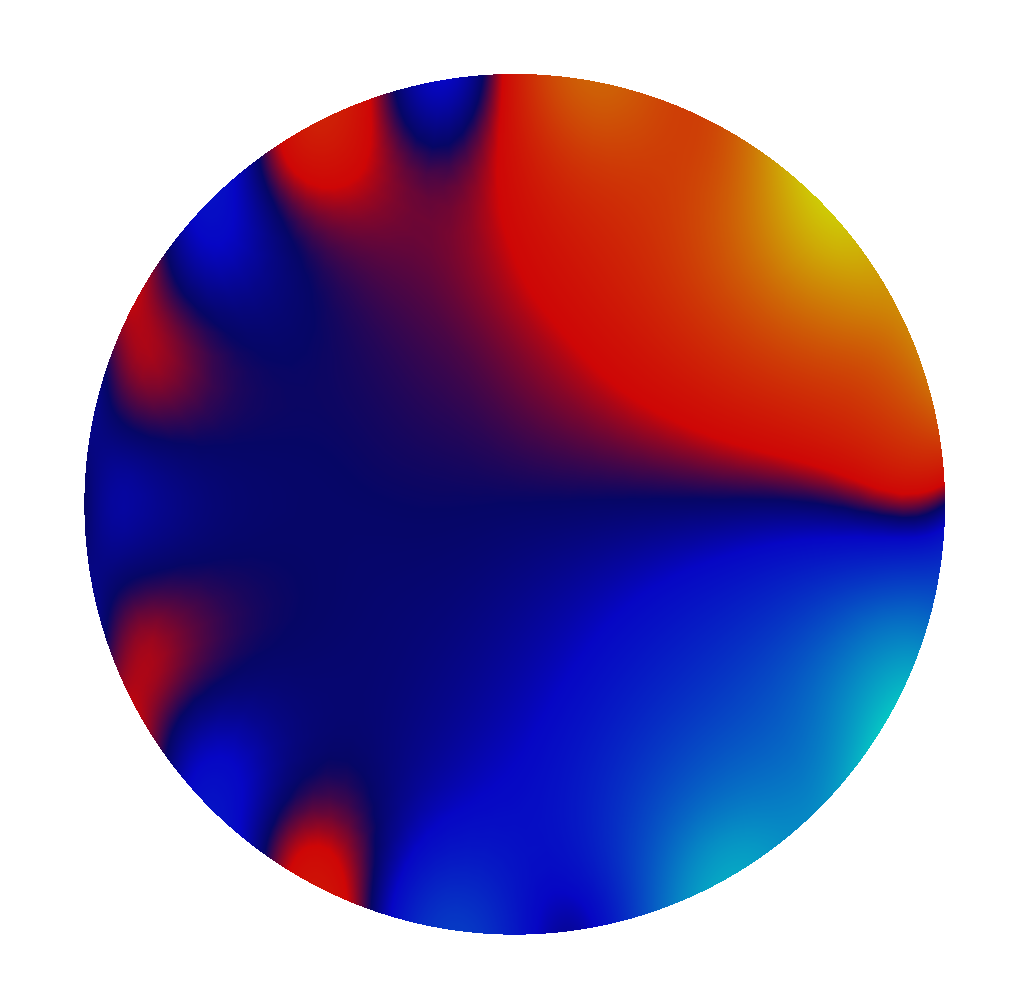} \\
    $0.25$ &
    \includegraphics[width=0.25\textwidth, valign=c]{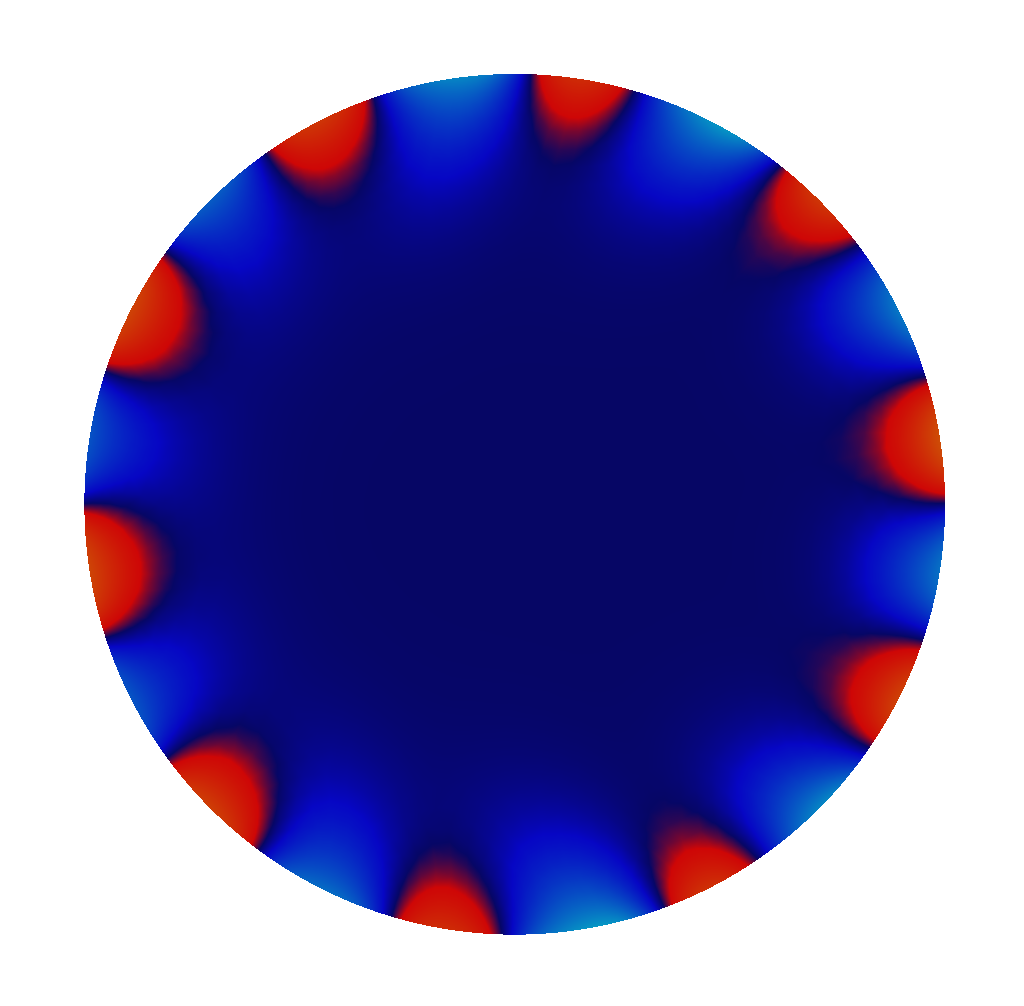} &
    \includegraphics[width=0.25\textwidth, valign=c]{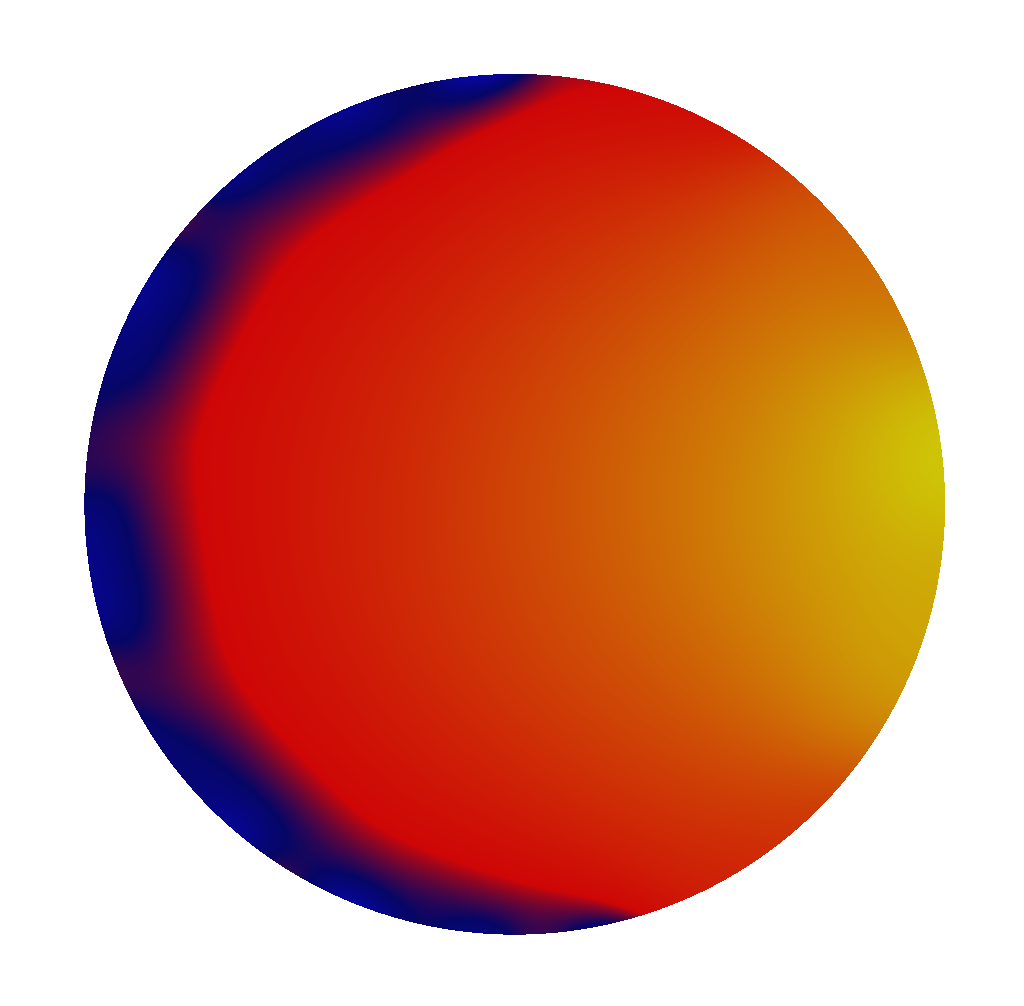} &
    \includegraphics[width=0.25\textwidth, valign=c]{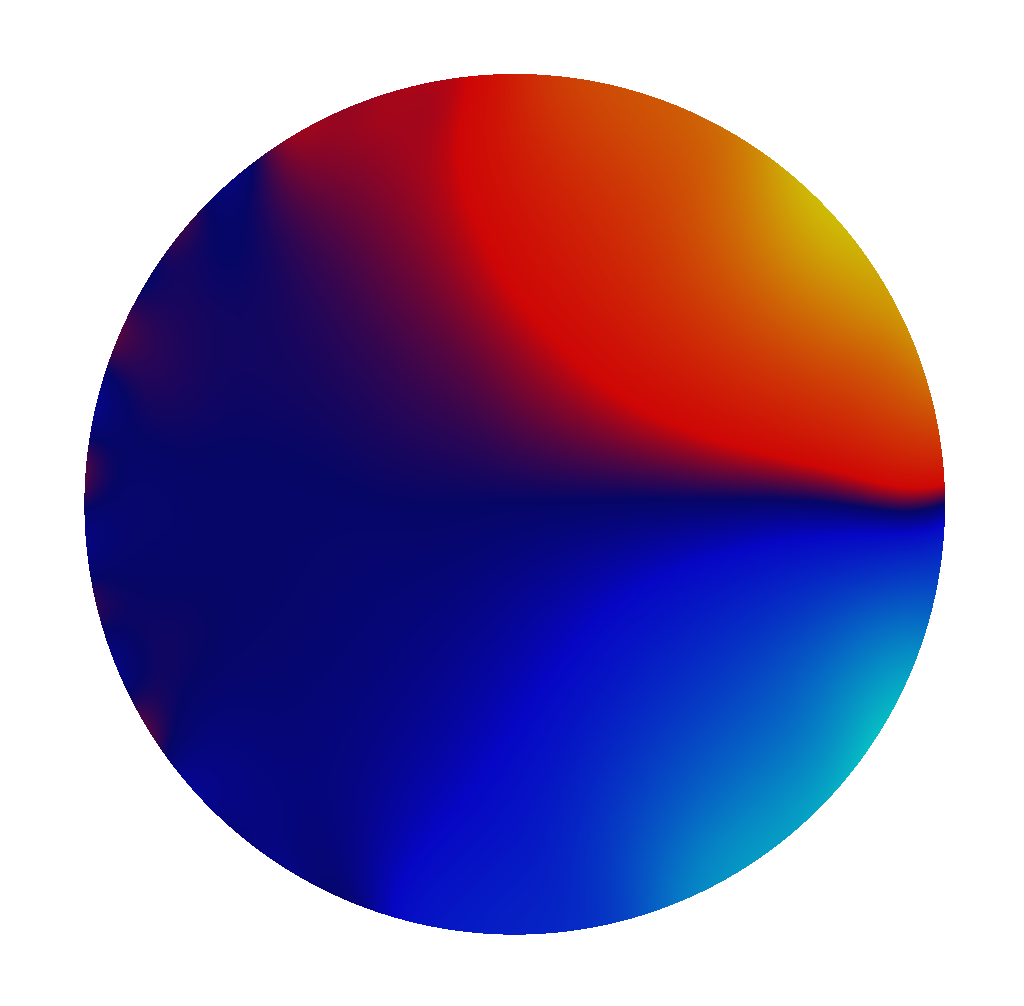} \\
    $0.50$ &
    \includegraphics[width=0.25\textwidth, valign=c]{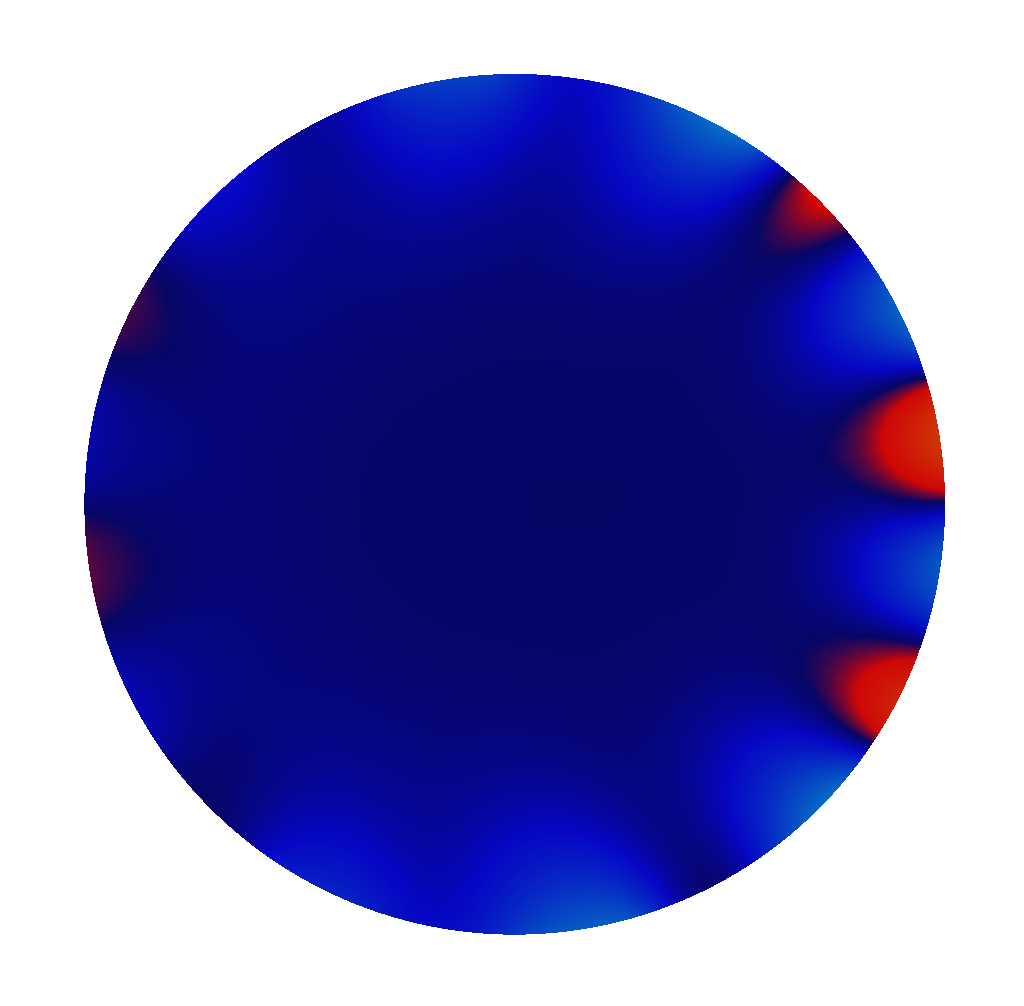} &
    \includegraphics[width=0.25\textwidth, valign=c]{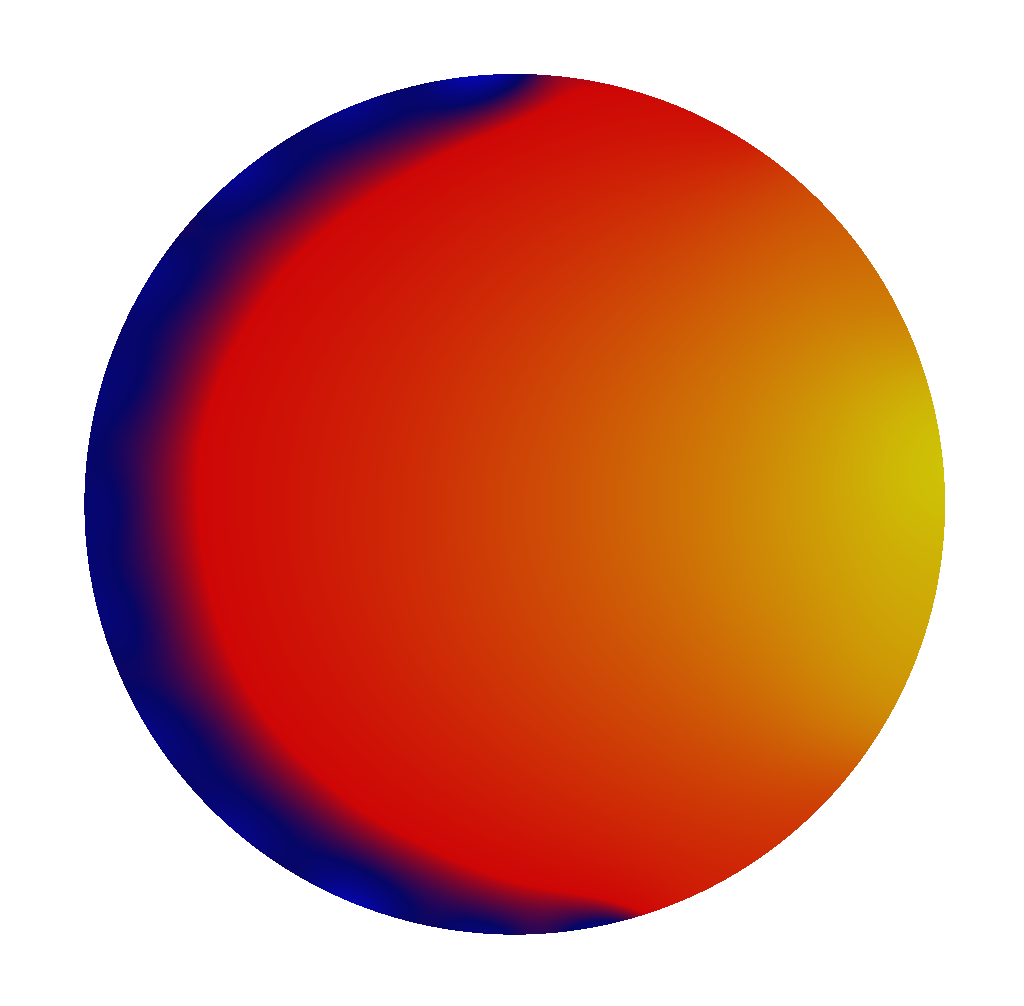} &
    \includegraphics[width=0.25\textwidth, valign=c]{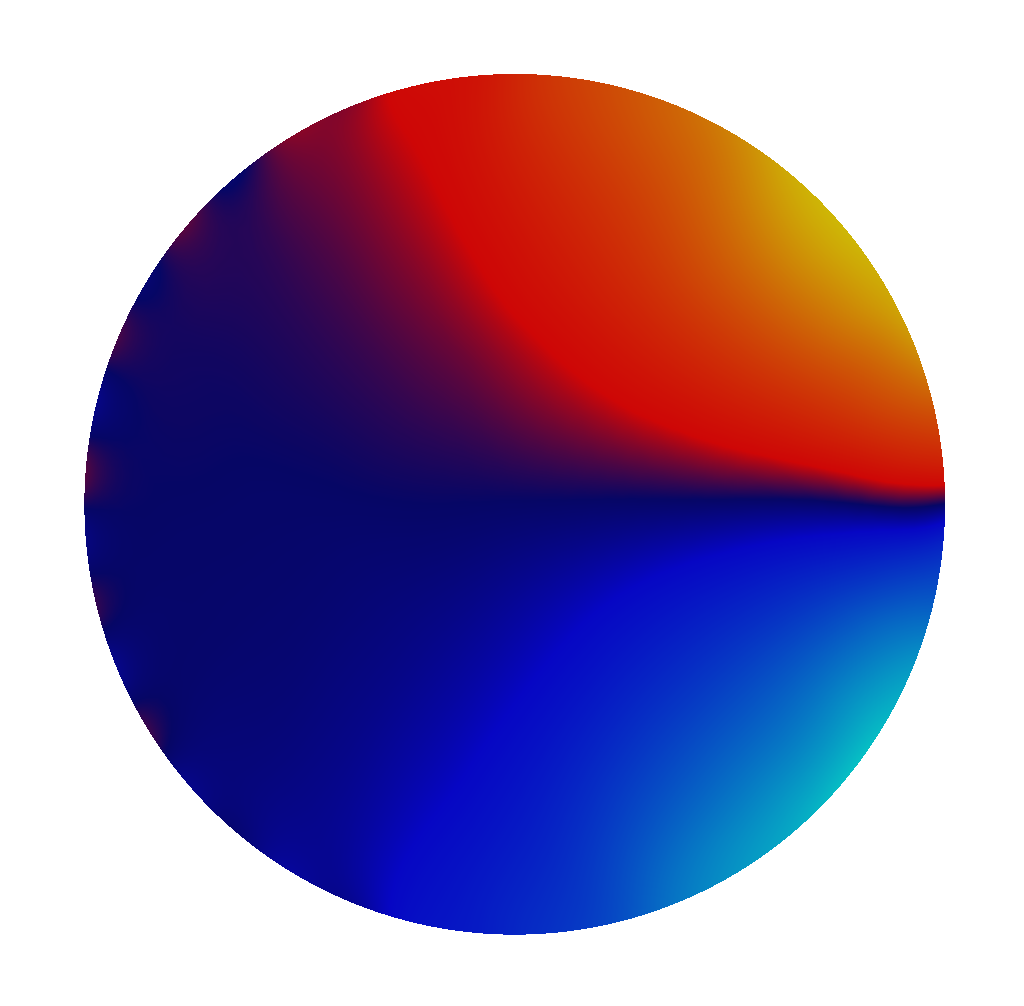} \\
    $0.75$ &
    \includegraphics[width=0.25\textwidth, valign=c]{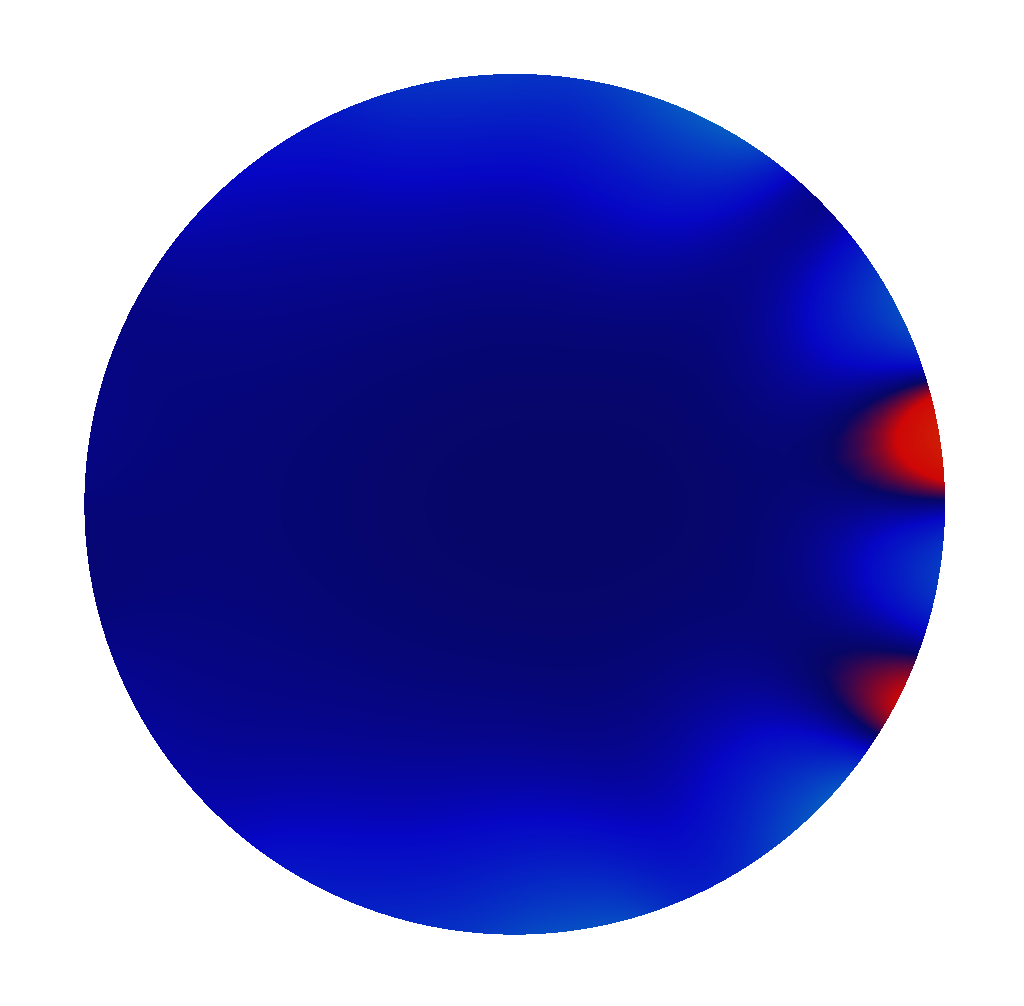} &
    \includegraphics[width=0.25\textwidth, valign=c]{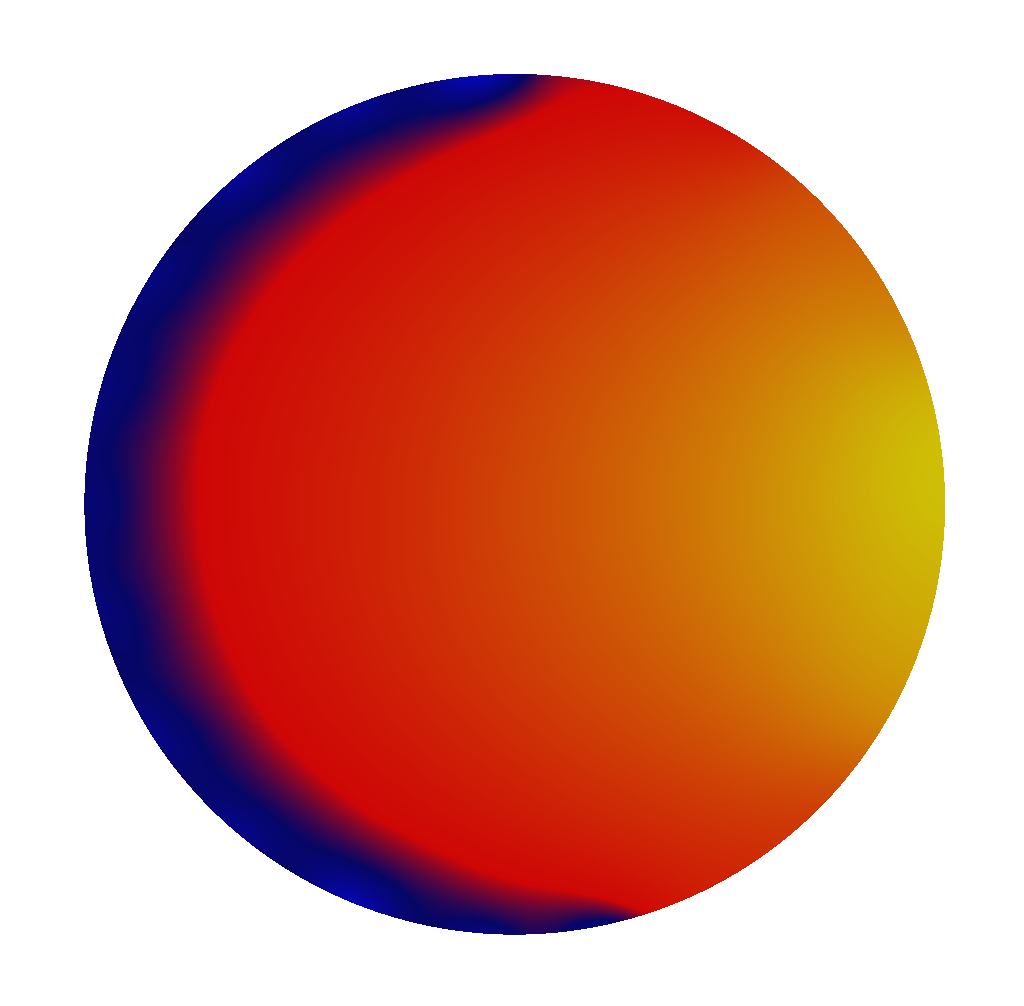} &
    \includegraphics[width=0.25\textwidth, valign=c]{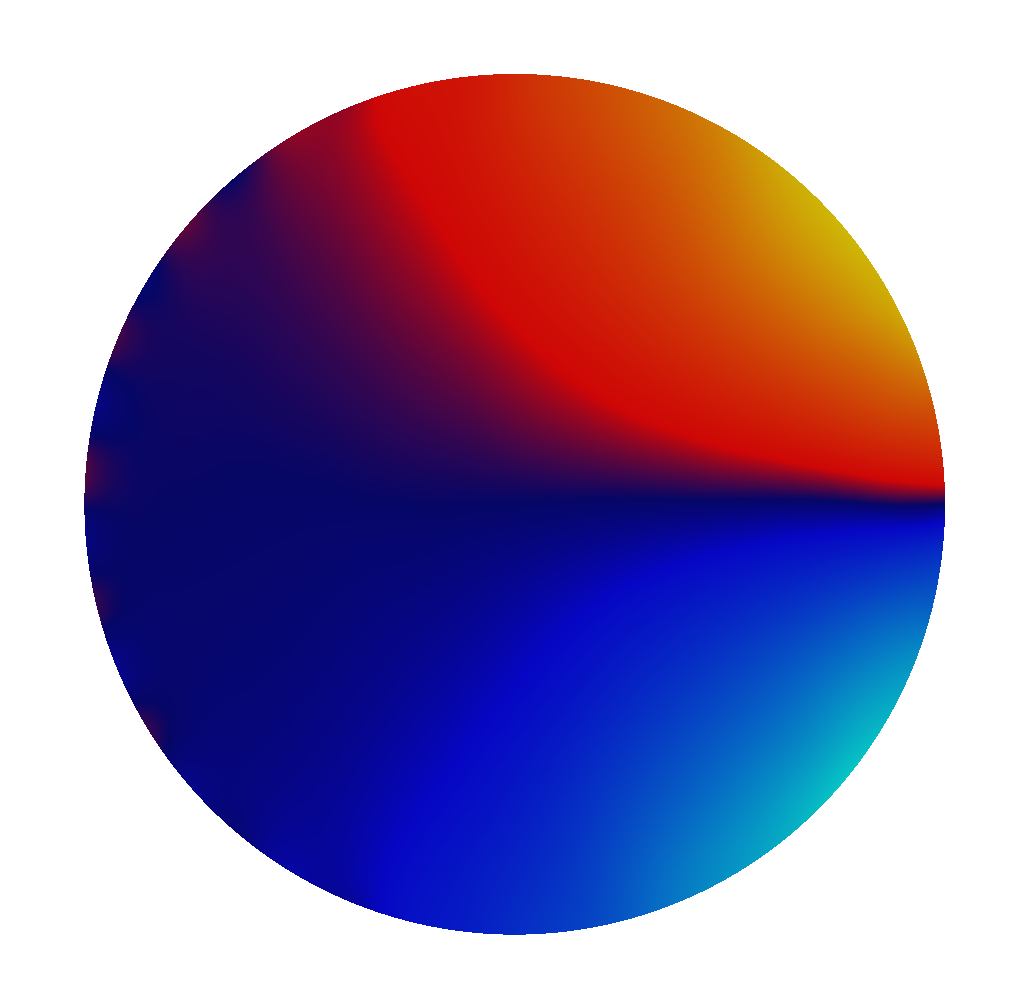} \\
    $1.00$ &
    \includegraphics[width=0.25\textwidth, valign=c]{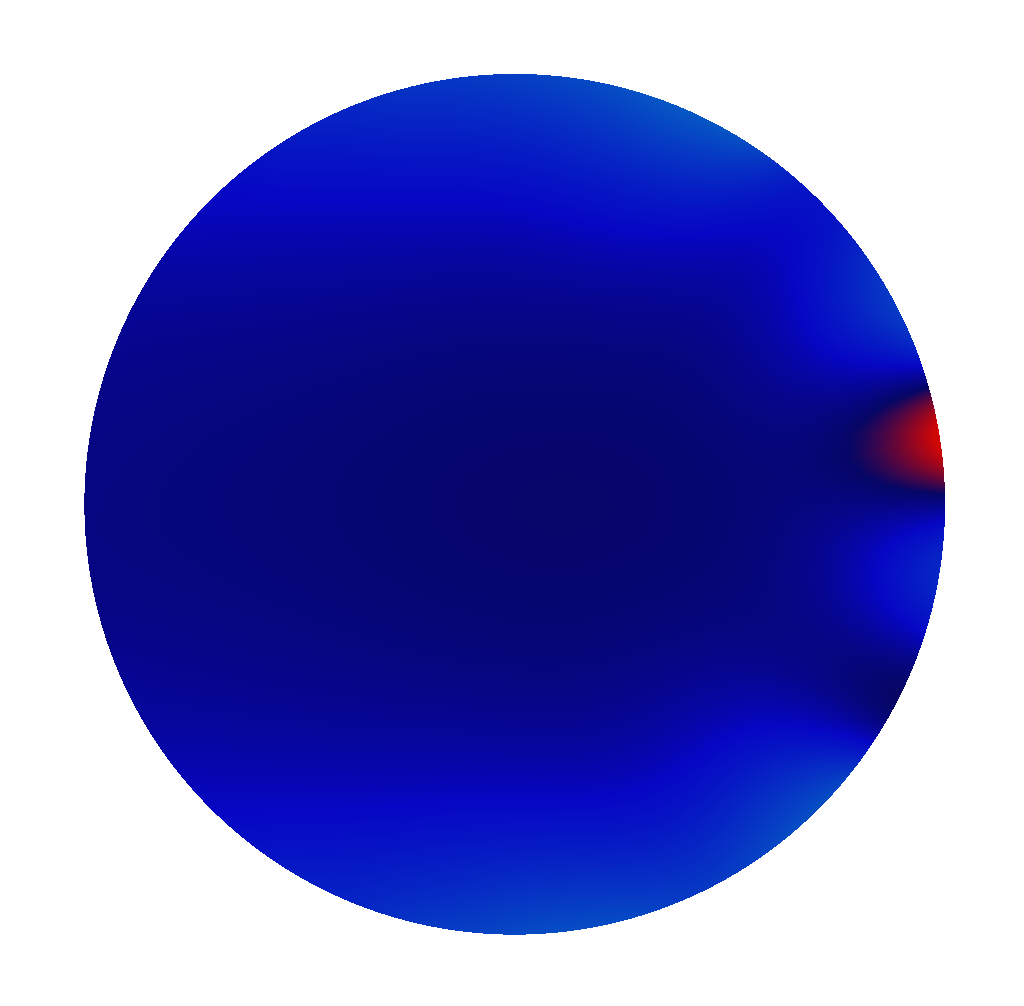} &
    \includegraphics[width=0.25\textwidth, valign=c]{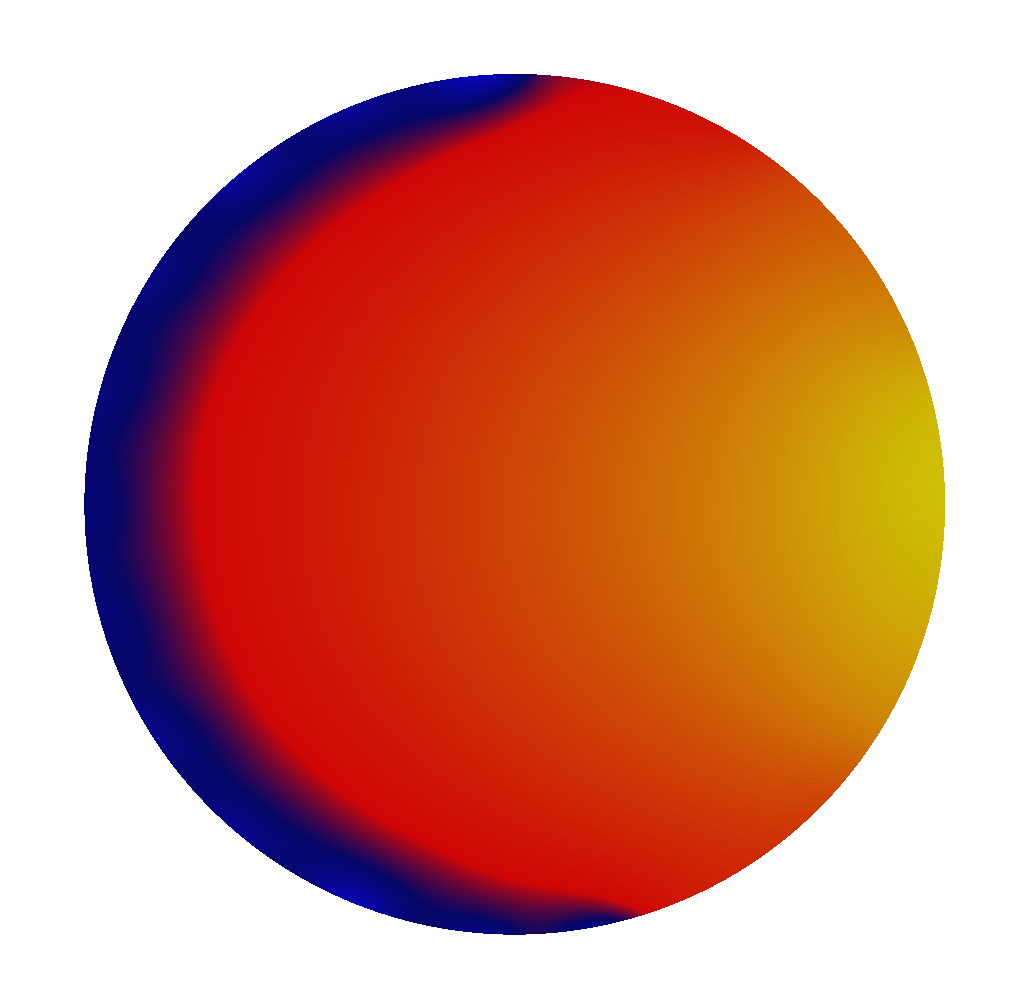} &
    \includegraphics[width=0.25\textwidth, valign=c]{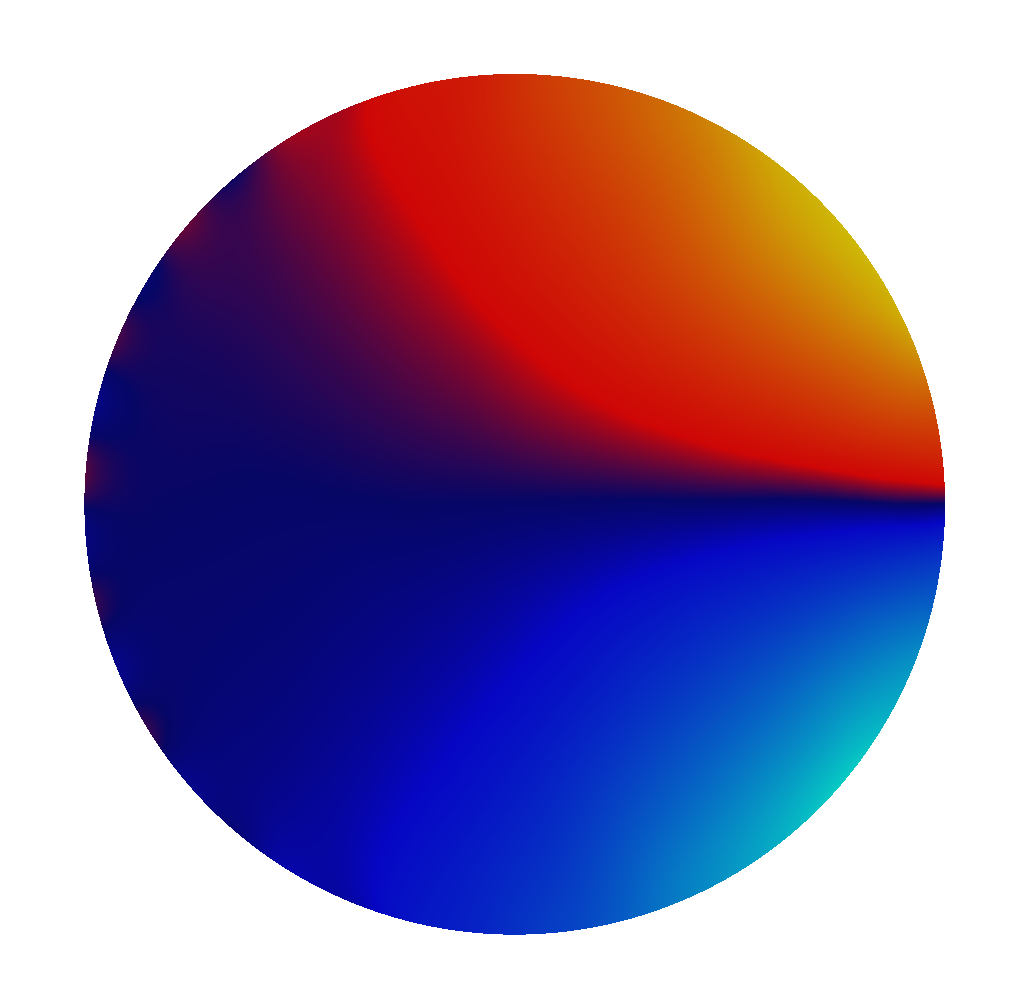} \\
     &
    \pgfplotscolorbardrawstandalone[
        colormap={osmosisc}{
            rgb(0.67374402284622203)=(0, 1, 1)
            rgb(0.95)=(0, 0, 1)
            rgb(1.0)=(0, 0, 0.50196078431400004)
            rgb(1.05)=(1, 0, 0)
            rgb(1.5269199609756501)=(1, 1, 0)
        },
        point meta min=0.67374402284622203,
        point meta max=1.5269199609756501,
        colorbar horizontal,
        colorbar style={xtick={0.67374402284622203,1.,1.5269199609756501}},
        parent axis width/.initial=0.20\textwidth
    ] &
    \pgfplotscolorbardrawstandalone[
        colormap={osmosisda}{
            rgb(-0.098934702575206798)=(0, 1, 1)
            rgb(-0.05)=(0, 0, 1)
            rgb(0.0)=(0, 0, 0.50196078431400004)
            rgb(0.05)=(1, 0, 0)
            rgb(1.07123994827271)=(1, 1, 0)
        },
        point meta min=-0.098934702575206798,
        point meta max=1.07123994827271,
        colorbar horizontal,
        colorbar style={xtick={-0.098934702575206798,0.,1.07123994827271},
                        xticklabels={{\hspace{-0.24cm}$-0.1$},{\hspace{0.2cm}$0$},{$1.07$}}},
        parent axis width/.initial=0.20\textwidth
    ] &
    \pgfplotscolorbardrawstandalone[
        colormap={osmosisdb}{
            rgb(-0.43490800261497498)=(0, 1, 1)
            rgb(-0.05)=(0, 0, 1)
            rgb(0.0)=(0, 0, 0.50196078431400004)
            rgb(0.05)=(1, 0, 0)
            rgb(0.45362100005149802)=(1, 1, 0)
        },
        point meta min=-0.43490800261497498,
        point meta max=0.45362100005149802,
        colorbar horizontal,
        colorbar style={xtick={-0.43490800261497498,0.,0.45362100005149802}},
        parent axis width/.initial=0.20\textwidth
    ] 
\end{tabular}
\caption{Solution of the ALE discretization of the example problem for parameters $\alpha = 0.1$, $\beta = 0.1$, $\delta_1 = 1$, $\delta_2 = 1$.
Depicted are the concentration field $\huh{n}$ as well as the scalar components of the deformation
field $\Psih{n} - \id$ for different times $t$.}\label{fig:solution}
\end{figure}

The solution trajectories of the solution fields $\huh{n}$ and $\Psih{n}$ are visualized in
Figure~\ref{fig:solution} for $\mu^* = (0.1, 0.1, 1, 1)$.
The corresponding reconstruction $\uh{n}$ on the deformed domain $\Omega_{h}^n$ is shown in
Figure~\ref{fig:reconstruction}.

The discrete model and its reduction was implemented using \texttt{NGSolve} \cite{ngsolve} and \texttt{pyMOR} \cite{MilkRaveEtAl2016}.
All computations were performed on a single core of dual-socket Intel Xeon E5-2698 compute server
with 256GB RAM.

\subsection{Lagrangian vs. Eulerian viewpoint} \label{sec:numex:lagvseul}

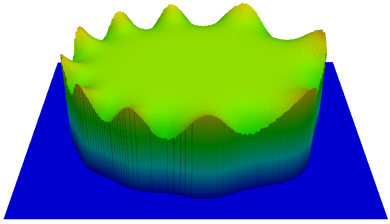
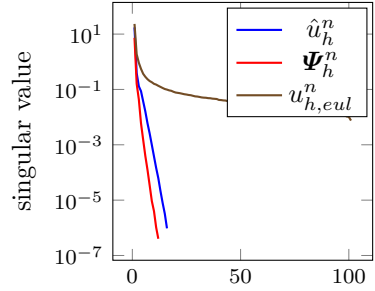
\begin{figure}[t] \centering
\begin{subfigure}[t]{0.45\textwidth}
\centering
\includegraphics[width=\textwidth]{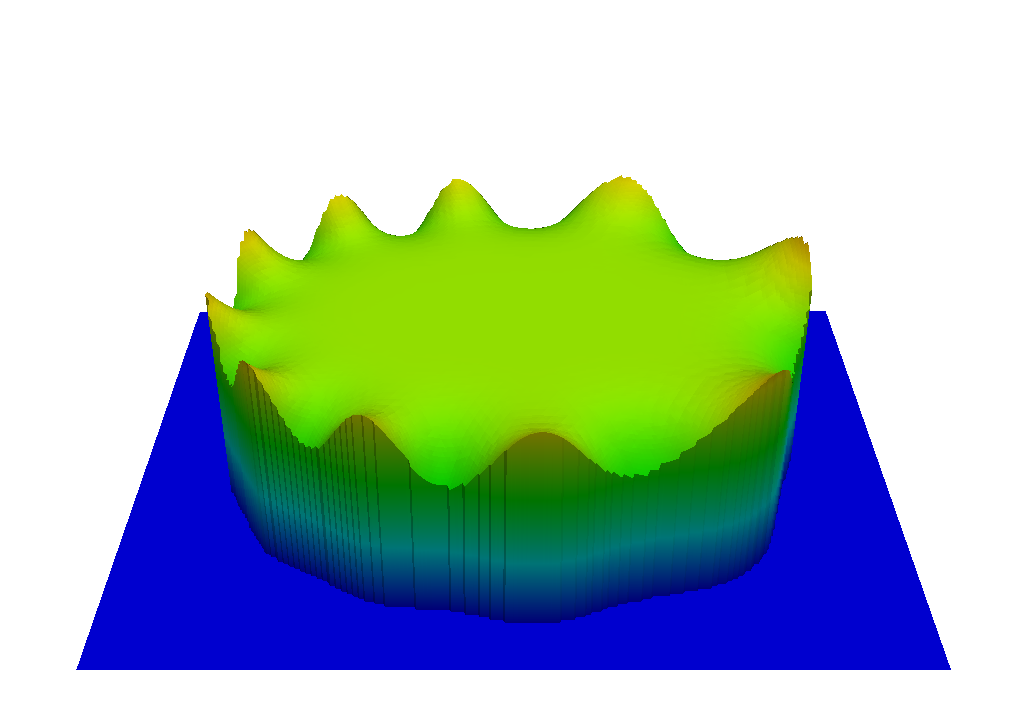}
\caption{Eulerian embedding of the concentration field trajectory $u_{h,eul}^n$ at time $t=0.25$.}\label{fig:eulerian}
\end{subfigure}
\qquad
\begin{subfigure}[t]{0.45\textwidth}
\begin{tikzpicture}
     \begin{semilogyaxis}[
             width=5.0cm,
             height=5.0cm,
             ylabel={singular value},
             ytick={1e1, 1e-1, 1e-3, 1e-5, 1e-7},
             tick label style={font=\footnotesize},
            ]
            \addplot+[mark=none,thick] table[x index=0, y index=1] {data/u_svals};
            \addplot+[mark=none,thick] table[x index=0, y index=1] {data/d_svals};
            \addplot+[mark=none,thick] table[x index=0, y index=1] {data/u_ref_svals};
            \legend{$\huh{n}$, $\Psih{n}$, $u_{h,eul}^n$}
     \end{semilogyaxis}
\end{tikzpicture}
\caption{Singular value decay of the concentration and transformation field trajectories $\huh{n}$, $\Psih{n}$ vs.~%
Eulerian concentration field trajectory $u_{h, eul}^n$.}\label{fig:svals}
\end{subfigure}
\caption{Eulerian embedding and singular value decay for the solution trajectory depicted in Figure
\ref{fig:solution}.}
\end{figure}

To assess the viability of our approach, consider for the same parameter $\mu^*$ the embeddings \begin{equation} \label{eq:lhn}
    \Lambda_h^n: L^2(\Omega^n_{h}) \longrightarrow L^2(\mathbb{R}^2)
\end{equation}
given by extending functions on $\Omega_h^n$ with zero outside of $\Omega_h^n$.
To numerically approximate this mapping we consider a fixed reference mesh on a sufficiently large
domain containing $\bigcup_{0\leq n \leq N} \Omega^n_h$ and compute the first order finite element
functions $u^n_{h, eul}$ on this mesh given by Lagrange interpolation of $\uh{n}$ (cf.\ Figure~\ref{fig:eulerian}).
These functions $u^n_{h, eul}$ can be seen as representatives for solution trajectories of
discretizations of \eqref{eq:osmosisproblem} that take an Eulerian viewpoint, as opposed to the presented
Lagrangian formulation, e.g.\ phase-field discretizations with the size of the diffuse interface tending
to zero.

In Figure~\ref{fig:svals} we compare the singular value decay of the Lagrangian solution field trajectories $\huh{n}$, $\Psih{n}$
to the singular value decay of the Eulerian concentration field trajectory $u^n_{h, eul}$.
As expected, due to the smooth time dependence of $\huh{n}$ and $\Psih{n}$, both trajectories show a
much faster singular value decay compared to the $u^n_{h, eul}$ trajectory, which is non-differentiable in
time due to its moving jump at the boundary of $\Omega_{h}^n$.

\begin{remark} \label{rem:unffemnohyperred}
We notice that exchanging the extensions $\Lambda_h^n$ in \eqref{eq:lhn}, e.g.\ applying recent ideas from unfitted finite elements \cite{LO_ARXIV_2018} to compute smooth extensions of the concentration field outside of $\Omega(t)$, may improve the approximability for the Eulerian formulation.
    However, even if the approximability can be improved drastically, it is unclear how hyper-reduction can be effectively applied (i.e.\ with a small number of interpolation points) in such an unfitted setting.

In \cite{BalajewiczFarhat2014} a different approach is taken, effectively considering the approximation
problem over multiple spaces, by truncating approximating vectors defined on
$L^2(\mathbb{R}^n)$ to the respective domain $\Omega(t)$.
\end{remark}

\subsection{Parametric model order reduction} \label{sec:numex:pmor}
\begin{figure}
\centering
\begin{subfigure}[5]{0.45\textwidth}
\begin{tikzpicture}
    \begin{axis}[
             grid=major,
             xlabel={dimension},
             ylabel={$\varepsilon_{rb}$},
             legend entries={$\Uru$, $\bUrPsi$, $\bUrG$},
             ymode=log,
             height=5cm,
             width=5cm,
            ]
        \addplot table[y=RBTOL, x=UDIM] {data/out.dat};
        \addplot table[y=RBTOL, x=DDIM] {data/out.dat};
        \addplot table[y=RBTOL, x=QNEDIM] {data/out.dat};
    \end{axis}
\end{tikzpicture}
\end{subfigure}
\quad
\begin{subfigure}[5]{0.45\textwidth}
\begin{tikzpicture}
    \begin{axis}[
             grid=major,
             xlabel={$M_1, \dots, M_7$},
             ylabel={$\varepsilon_{ei}$},
             legend entries={$c_1$, $c_2$, $c_3$, $c_4$, $c_5$, $c_6$, $c_7$},
             ymode=log,
             height=5cm,
             width=5cm,
             empty line=none,
            ]
        \addplot table[y=EITOL, x=MQ1DIM] {data/out.dat};
        \addplot table[y=EITOL, x=MQ2DIM] {data/out.dat};
        \addplot table[y=EITOL, x=MUDIM] {data/out.dat};
        \addplot table[y=EITOL, x=AU1DIM] {data/out.dat};
        \addplot table[y=EITOL, x=AU2DIM] {data/out.dat};
        \addplot table[y=EITOL, x=LQ1DIM] {data/out.dat};
        \addplot table[y=EITOL, x=LQ2DIM] {data/out.dat};
    \end{axis}
\end{tikzpicture}
\end{subfigure}
\caption{Reduced space dimensions (left) and number of empirical interpolation points (right) vs.\ training
tolerance for the numerical experiment.}\label{fig:basis_sizes}
\end{figure}
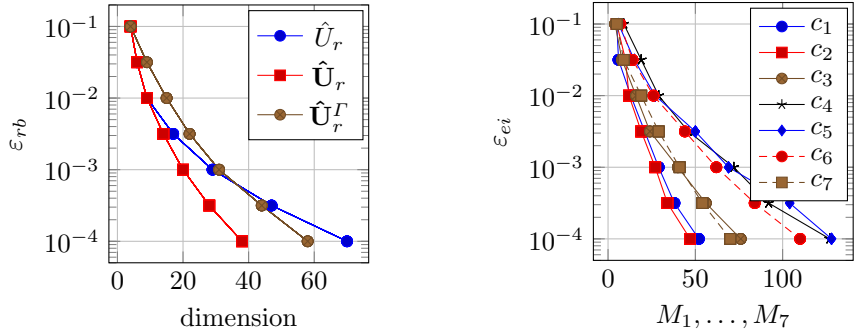

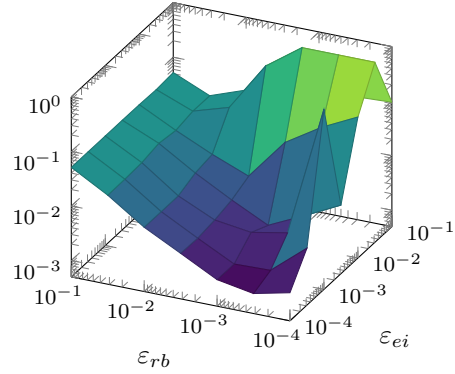
\begin{figure}[t]
\centering
\begin{subfigure}[t]{0.48\textwidth}\centering
 \begin{tikzpicture}
     \begin{axis}[
             title={Model Order Reduction Error -- $\hurmu{n}$},
             legend style={legend pos=north east,anchor=north west,xshift=9},
             width=5.8cm,
             height=5.8cm,
             colormap/viridis,
             x dir=reverse,
             xmode=log,
             ymode=log,
             zmode=log,
             zmin=5e-4,
             zmax = 1.,
             xlabel={$\varepsilon_{rb}$},
             ylabel={$\varepsilon_{ei}$},
             tick label style={font=\footnotesize},
             xticklabel style={xshift=-1, yshift=2},
             yticklabel style={xshift=6, yshift=9},
             xlabel shift=-3,
             ylabel shift=-3,
            ]
            \addplot3[surf] table[x=RBTOL, y=EITOL, z expr={min(\thisrow{UERROR},1)}] {data/out.dat};
     \end{axis}
\end{tikzpicture}
 \end{subfigure}
\quad
\begin{subfigure}[t]{0.48\textwidth}\centering
 \begin{tikzpicture}
     \begin{axis}[
             title={Model Order Reduction Error -- $\Psirmu{n}$},
             legend style={legend pos=north east,anchor=north west,xshift=9},
             width=5.8cm,
             height=5.8cm,
             colormap/viridis,
             x dir=reverse,
             xmode=log,
             ymode=log,
             zmode=log,
             zmax = 1.,
             zmin=5e-4,
             xlabel={$\varepsilon_{rb}$},
             ylabel={$\varepsilon_{ei}$},
             tick label style={font=\footnotesize},
             xticklabel style={xshift=-1, yshift=2},
             yticklabel style={xshift=6, yshift=9},
             xlabel shift=-3,
             ylabel shift=-3,
            ]
            \addplot3[surf] table[x=RBTOL, y=EITOL, z expr={min(\thisrow{DERROR},1)}] {data/out.dat};
     \end{axis}
\end{tikzpicture}
\end{subfigure}
\caption{Model order reduction errors for the numerical experiment. Depicted is the maximum relative error over all
$0 \leq n \leq N$ and $\mu \in \mathcal{S}_{test}$ for a test set $\mathcal{S}_{test} \subset \mathcal{P}$ of
100 randomly chosen parameters. The errors were truncated at a maximum value of 1 to improve the readability of the plots.}\label{fig:mor_errors}
\end{figure}

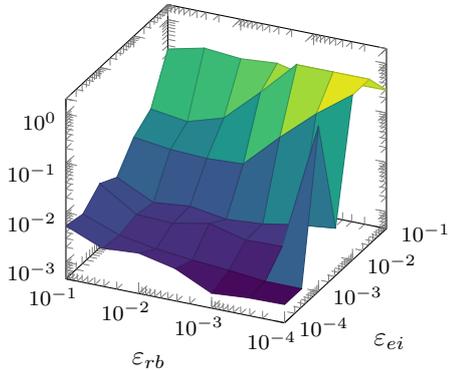
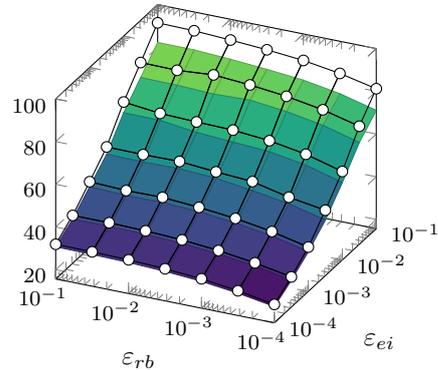
\begin{figure}[t]
\centering
\begin{subfigure}[t]{0.48\textwidth}\centering
 \begin{tikzpicture}
     \begin{axis}[
             title={Rel. Mass Conservation Error},
             legend style={legend pos=north east,anchor=north west,xshift=9},
             width=5.8cm,
             height=5.8cm,
             colormap/viridis,
             x dir=reverse,
             xmode=log,
             ymode=log,
             zmode=log,
             xlabel={$\varepsilon_{rb}$},
             ylabel={$\varepsilon_{ei}$},
             tick label style={font=\footnotesize},
             xticklabel style={xshift=-1, yshift=2},
             yticklabel style={xshift=6, yshift=9},
             xlabel shift=-3,
             ylabel shift=-3,
            ]
            \addplot3[surf] table[x=RBTOL, y=EITOL, z expr={min(\thisrow{TCERRORNC},1)}] {data/out.dat};
     \end{axis}
\end{tikzpicture}
    \caption{Maximum relative errors in mass conservation \eqref{eq:cons_error} for $\mu \in \mathcal{S}_{test}$ when using
    \eqref{eq:rom_u_ei} for the concentration update. The errors were truncated at a maximum value of 1 to improve the readability of the plots.}\label{fig:conserror}
\end{subfigure}
\quad
\begin{subfigure}[t]{0.48\textwidth}\centering
 \begin{tikzpicture}
     \begin{axis}[
             title={Model Order Reduction Speedup},
             width=5.8cm,
             height=5.8cm,
             colormap/viridis,
             x dir=reverse,
             ztick distance=20,
             zmax=100,
             xmode=log,
             ymode=log,
             xlabel={$\varepsilon_{rb}$},
             ylabel={$\varepsilon_{ei}$},
             tick label style={font=\footnotesize},
             xticklabel style={xshift=-1, yshift=2},
             yticklabel style={xshift=6, yshift=9},
             xlabel shift=-3,
             ylabel shift=-3,
            ]
            \addplot3[surf] table[x=RBTOL, y=EITOL, z=SU] {data/out.dat};
            \addplot3[mesh,mark=*,draw=black,fill=white] table[x=RBTOL, y=EITOL, z=SUNC] {data/out.dat};
     \end{axis}
\end{tikzpicture}
\caption{Median of model order reduction speedup over $\mu \in \mathcal{S}_{test}$. Surface plot: concentration update
using \eqref{eq:rom_u_ei_cons}; mesh plot: concentration update using \eqref{eq:rom_u_ei}.}\label{fig:speedup}
\end{subfigure}
\caption{Total mass conservation errors and model order reduction speedups for the numerical experiment (cf.\
Figure~\ref{fig:mor_errors}).}
\end{figure}

To test the parameterized model order reduction approach discussed in Section~\ref{sec:MOR},
we consider the parameter domain
\begin{equation}
    (\alpha, \beta, \delta_1, \delta_2) \in \mathcal{P} := [0.1, 1] \times [0.001, 0.1] \times [0, 1]^2 \subset
    \mathbb{R}^4,
\end{equation}
of which we choose a training set $\mathcal{S}_{train} \subset \mathcal{P}$ of $3^4$ equidistant parameters.
From the corresponding snapshot trajectories we compute the reduced approximation spaces $\Uru$, $\bUrPsi$, $\bUrG$
and empirical interpolations of $c_1, \ldots, c_7$ for varying relative training error tolerances
$\varepsilon_{rb}$, $\varepsilon_{ei}$ (cf.\ Figure~\ref{fig:basis_sizes}).
To assess the quality of the resulting ROMs, we compute the maximum relative model order reduction
errors for 100 randomly chosen test parameters $\mathcal{S}_{test} \subset \mathcal{P}$ (cf.\ Figure~\ref{fig:mor_errors}). 
We can observe an exponential error decay in both solution variables for simultaneously decreasing
training error tolerances.
For $\varepsilon_{rb} = \varepsilon_{ei} = 10^{-3}$ we observe errors of $1.36 \cdot 10^{-2}$ for $\hurmu{n}$ and of $2.99 \cdot 10^{-3}$
for $\Psirmu{n}$.
For $\varepsilon_{ei} \gg \varepsilon_{rb}$ we observe the usual instability of ROMs employing
empirical interpolation for hyper-reduction.
The computational speedup of the ROM over the finite element ALE discretization is shown in
Figure~\ref{fig:speedup} (surface plot). For $\varepsilon_{rb} = \varepsilon_{ei} = 10^{-3}$ the speedup is 41.

Next we consider the effect of not enforcing total mass conservation of the ROM:
Whereas all previous computations were performed using the mass conservative concentration update equation
\eqref{eq:rom_u_ei_cons}, we now compute solutions of the ROM obtained using \eqref{eq:rom_u_ei}
for the concentration update and consider the maximum relative mass conservation error
\begin{equation}\label{eq:cons_error}
    \frac{|\int_{\Omega_{r,\mu}^N} \hurmu{n} d\bx - \int_{\Omega_{r,\mu}^0} \hurmu{0} d\bx|}
         {\int_{\Omega_{r,\mu}^0} \hurmu{0} d\bx}
\end{equation}
over all $\mu \in \mathcal{S}_{test}$ (Figure~\ref{fig:conserror}).
We observe that the error decays with, but is almost always larger than $\varepsilon_{ei}$, whereas
it is mostly independent of $\varepsilon_{rb}$.
At the same time, computing \eqref{eq:rom_u_ei} instead of \eqref{eq:rom_u_ei_cons} is only slightly
faster (Figure~\ref{fig:speedup}, mesh plot).
For $\varepsilon_{rb} = \varepsilon_{ei} = 10^{-3}$ the speedup grows from 41 to 43.
The model order reduction errors are of the same order of magnitude.

\begin{figure}[t] \centering
 \begin{tikzpicture}
     \begin{groupplot}[
             group style={group size=3 by 1,
                          horizontal sep=1.2cm,
                          vertical sep=1.5cm,},
             width=4.4cm,
             height=4.4cm,
             colormap/viridis,
             zmin=0,
             zmax=3e-3,
             /tikz/font=\small,
             xlabel={$\delta_1$},
             ylabel={$\delta_2$},
     ]
         \nextgroupplot[title={$t = 0.25$}]
            \addplot3[surf, mesh/rows=50] table[x=INIT0, y=INIT1, z=V1] {data/variance.dat};
         \nextgroupplot[title={$t = 0.5$}]
            \addplot3[surf, mesh/rows=50] table[x=INIT0, y=INIT1, z=V2] {data/variance.dat};
         \nextgroupplot[title={$t = 0.75$}]
            \addplot3[surf, mesh/rows=50] table[x=INIT0, y=INIT1, z=V3] {data/variance.dat};
     \end{groupplot}
\end{tikzpicture}
\caption{Variance $V^n_{r,\mu}$ \eqref{eq:variance} of the concentration field $\hurmu{n}$ for $50 \times 50$
combinations of the deformation parameters $\delta_1$, $\delta_2$, $\alpha = 0.1$, $\beta = 0.1$.}\label{fig:variance}
\end{figure}
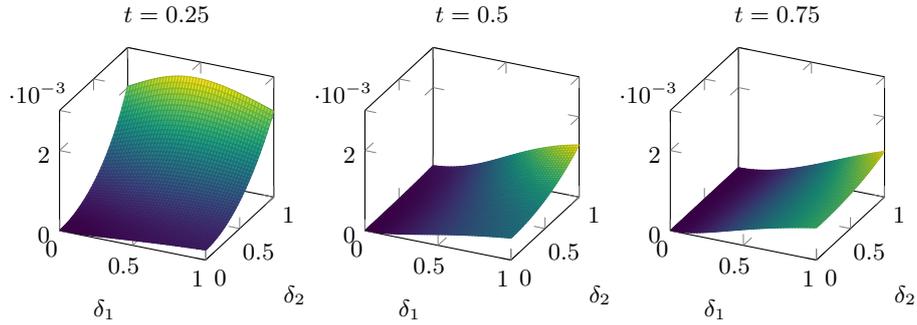

One important use case of parametric model order reduction is to quickly compute a certain
quantity of interest from the state-space solution of the ROM.
As an example, we here consider the variance of the concentration field $\hurmu{n}$
\begin{equation}\label{eq:variance}
    V_{r,\mu}^n := \frac{\int_{\Omega_{r,\mu}^n} (\hurmu{n} - \bar{u}^n_{r,\mu})^2 d\bx}
                        {\int_{\Omega_{r,\mu}^n} 1 d\bx}, \qquad
    \bar{u}_{r,\mu}^n := \frac{\int_{\Omega_{r,\mu}^n} u^n_{r,\mu} d\bx}
                              {\int_{\Omega_{r,\mu}^n} 1 d\bx},
\end{equation}
which can be computed from the ROM as
\begin{equation}\label{eq:variance_comp}
    V_{r,\mu}^n = \frac{\tilde{a}_3(\hurmu{n} - \bar{u}_{r, \mu}^n, \hurmu{n} - \bar{u}_{r, \mu}^n; \Psirmu{n})}
                        {\tilde{a}_3(1, 1; \Psirmu{n})}, \qquad
    \bar{u}_{r,\mu}^n = \frac{\tilde{a}_3(\hurmu{n}, 1; \Psirmu{n})}
                              {\tilde{a}_3(1, 1; \Psirmu{n})}.
\end{equation}
In Figure~\ref{fig:variance} we have computed $V_{r, \mu}$ for $\varepsilon_{rb} = \varepsilon_{ei} = 10^{-3}$
and $50 \times 50$ equidistant values of $\delta_1$ and $\delta_2$ whereas $\alpha = \beta = 0.1$ were fixed.
The computation of these $2,500$ outputs took 36 minutes.

\section*{Conclusion and outlook} \label{sec:conclusion}
We presented the application of projection-based model order reduction for a model free boundary problem with unknown
evolving geometry.
By using ALE mappings to describe the geometry evolution of the model, we were able to apply standard reduced basis
and empirical interpolation techniques.
An appropriate time discretization allowed us to obtain a globally mass conservative Galerkin ROM 
by including the constant functions in the reduced state space.
Through a rank-one modification of the empirically interpolated reduced mass matrix we could preserve this property
in the final fully online-efficient Galerkin-EI ROM.

In this work we have focussed on the reduced order modeling aspects specific to free boundary problems.
The integration of more advanced techniques such as greedy basis generation based on a posteriori
error indicators for the ROM should be straightforward.

We believe that the methodology can also be applied to different and even more complex free boundary problems. 
However, an obvious limitation of the presented approach lies in the dependency on one reference domain w.r.t.\ which all other domains can be expressed. For problems with large deformations or topology changes this will be insufficient and different discretizations with suitable reduced order modeling have to be considered.

Without changing the general approach, remeshing in the high-dimensional problem could be considered
by considering multiple meshes on the reference domain $\hat\Omega$ and using reduced basis techniques for
mesh-adaptive schemes such as \cite{Yano2014,UllmannRotkvicEtAl2016,AliSteihEtAl2014,GraessleHinze2017}.
However in the presence of large deformations, $\Omega(t)$-dependent norms will have to be considered in the
schemes in contrast to the currently fixed function spaces norms w.r.t.\ the reference domain.

Very attractive for the discretization would be the use of unfitted (or embedded) discretizations. In such a setting a
naive approximation fails as we have seen in Section \ref{sec:numex:lagvseul}. It would be
interesting if the slow decay of the Kolmogorov $n$-widths could be repaired by choosing a suitable embedding
to define a joint approximation space (see also Remark \ref{rem:unffemnohyperred}).
How hyper-reduction could be effectively applied in such an approach is unclear, however.

\section*{Software availability}
The source code used to produce the numerical results in Section~\ref{sec:numex} can be obtained
from \url{https://doi.org/10.5281/zenodo.1232550} under an open source license.

\end{document}